\documentclass[12pt, psamsfonts]{amsart}
\usepackage{amsmath}
\usepackage{amsthm}
\usepackage{amssymb}
\usepackage{amscd}
\usepackage{amsfonts}
\usepackage{amsbsy}
\usepackage{epsfig}

\theoremstyle{plain}
\newtheorem{thm}{Theorem}[section]

\newtheorem{lem}[thm]{Lemma}

\newtheorem{prop}[thm]{Proposition}
\newtheorem{claim}{Claim}
\newtheorem*{claim*}{Claim}

\theoremstyle{definition}

\newtheorem{rmk}[thm]{Remark}
\newcommand{\st}{such that }

\newcommand{\es}{\emptyset }
\newcommand{\bd}{\partial }

\newcommand{\sm}{\setminus}
\newcommand{\OOO}{\emptyset}

\newcommand{\RRR}{\mathbb{R}}
\newcommand{\CCC}{\mathbb{C}}
\newcommand{\NNN}{\mathbb{N}}

\renewcommand{\le}{\leqslant}
\renewcommand{\ge}{\geqslant}

\begin{document}

\title{Distortion bounds for $C^{2+\eta}$ unimodal maps}
\author{Mike Todd}


\begin{abstract}
We obtain estimates for derivative and cross--ratio distortion for
$C^{2+\eta}$ (any $\eta>0$) unimodal maps with non--flat critical
points.  We do not require any `Schwarzian--like' condition.

For two intervals $J \subset T$, the cross--ratio is defined as
the value
$$B(T,J):=\frac{|T||J|}{|L||R|}$$ where $L,R$ are the left and
right connected components of $T\setminus J$ respectively.  For an
interval map $g$ \st $g_T:T \to \RRR$ is a diffeomorphism, we
consider the cross--ratio distortion to be
$$B(g,T, J):=\frac{B(g(T),g(J))}{B(T,J)}.$$

We prove that for all $0<K<1$ there exists some interval $I_0$
around the critical point \st for any intervals $J \subset T$, if
$f^n|_T$ is a diffeomorphism and $f^n(T) \subset I_0$ then
$$B(f^n, T, J)> K.$$
Then the distortion of derivatives of $f^n|_J$ can be estimated
with the Koebe Lemma in terms of $K$ and $B(f^n(T),f^n(J))$.  This
tool is commonly used to study topological, geometric and ergodic
properties of $f$. This extends a result of Kozlovski.
\end{abstract}

\maketitle

\section{Introduction}
\label{sec:introcr}

In order to understand the long term behaviour of smooth dynamical
system $f:X \to X$ we must consider iterates of the map.  It is
useful to know how differently high iterates of the map $f^n$ act
on nearby points. For example we can try to estimate how wild the
derivative of iterates of the map is: we can consider the
distortion $\frac{Df^n(x)}{Df^n(y)}$ for $x,y$ in some small
interval $J$ where $f^n|_J$ is a diffeomorphism.  For one
dimensional maps, the Koebe Lemma is a tool we use to estimate
this.  Notice that this distortion can be rather wild when $f$ has
critical points.

An important condition we must assume in order to apply the Koebe
Lemma is that the map $f^n$ must increase cross--ratios. The type
of cross--ratio we use most is defined as follows.  For two
intervals $J \subset T$, the \emph{cross--ratio} is defined as the
value
$$B(T,J):=\frac{|T||J|}{|L||R|}$$ where $L,R$ are the left and
right connected components of $T\setminus J$ respectively.  For an
interval map $g$ \st $g_T:T \to \RRR$ is a diffeomorphism, the
main measure of cross--ratio distortion we use is given by
$$B(g,T, J):=\frac{B(g(T),g(J))}{B(T,J)}.$$
If we know that $B(f^n,T^*,J^*)\ge K>0$ for any $J^*\subset
T^*\subset T$ then we have uniform bounds on
$\frac{Df^n(x)}{Df^n(y)}$ for $x,y\in J$ depending on $K$ and
$B(f^n(T), f^n(J))$.  So we are able to estimate the distortion of
the derivative of $f^n$ using information on the distortion of the
cross--ratios.

A classical way of gaining information about the dynamics of an
interval map $f:[0,1] \to [0,1]$ with a critical point, is to take
a first return map to some well chosen interval $I$. If this map
has some diffeomorphic branches, we can estimate how well or how
badly the derivatives behave on branches using the Koebe Lemma as
above. This method is often used to give information on the
geometry and topology of the map and its iterates, see \cite{ms}.
This type of approach is also applied when considering the ergodic
properties of one dimensional maps. Often instead of first return
maps, certain \emph{inducing schemes} are applied in these cases,
see \cite{ms}.  The Koebe Lemma allows us to show that the
inducing schemes are expansive, and the Folklore Theorem can then
be used to derive ergodic absolutely continuous $f$--invariant
measures.

In order to apply the Koebe Lemma to $f^n|_T$ we need a lower
bound on cross ratio distortion of $f^n|_T$. In fact, a lower
bound $K=1$ is obtained whenever $f$ is $C^3$ and has negative
\emph{Schwarzian derivative}: that is
$$Sf:=\frac{D^3f}{Df}-\frac 32\left(\frac{D^2f}{Df}\right)^2$$
is negative wherever it is well defined. For applications it is
not so important that $f$ have negative Schwarzian, just that some
iterate of $f$ has negative Schwarzian on some small intervals.
Kozlovski showed \cite{kozsch} that for any $C^3$ unimodal map
with non--flat critical point (see the next section), if $I$ is a
small enough neighbourhood of the critical point and $f^n(x)\in I$
then $Sf^{n+1}(x)<0$. Therefore, for most practical purposes, for
example where first return maps or inducing schemes are used to
gain information about the dynamics, it is unnecessary to find the
sign of the Schwarzian derivative as long as the critical point is
non--flat.  Moreover, this result allowed Kozlovski to prove the
following.
\begin{thm}Suppose that $f$ is a $C^3$ unimodal map with non--flat
critical point whose iterates do not converge to a periodic
attractor. Then for any $0<K<1$, there is an interval $V$ around
the critical point \st if, for an interval $T$ and some $n>0$,
\begin{itemize}
\item $f^n|_T$ is monotone; and
\item each interval from the orbit $\{T, f(T), \ldots, f^n(T)\}$ is
contained in the domain of the first entry map to $V$,
\end{itemize}
then $$B(f^n, T, J) > K$$ where $J$ is any subinterval of $T$.
\label{thm:koz}
\end{thm}
This means that the Koebe Lemma can be applied to $f^n$ to get
estimates on the distortion of derivatives which only depend on
$B(f^n(T), f^n(J))$ (for first return maps or induced maps this
quantity is bounded whenever the branches have a `uniform
extension').  These results were extended to $C^3$ multimodal maps
with non--flat critical points in \cite{svarg}. Also, for $C^3$
unimodal maps with non--flat critical point, it is shown in
\cite{gss2} that an analytic coordinate change can create a map
which has first return maps with negative Schwarzian.

So how necessary is the negative Schwarzian condition to prove
dynamical results in `reasonable' cases?  Certainly it is useful
in determining the type of parabolic periodic points or bounding
the number of attracting cycles, see \cite{singer, ms}.  A natural
question to ask, and the one we consider in this paper, is: what
happens for unimodal maps with non--flat critical points which are
not $C^3$? Certainly the usual negative Schwarzian condition is no
use since it is not even defined. (Note that there is a
`Schwarzian--like' condition for $C^1$ maps -equivalent to the
negative Schwarzian condition when the map is $C^3$- see
\cite{preston1,ms}, but that need not hold in our case either.) We
show that Theorem~\ref{thm:koz} extends to the case of
$C^{2+\eta}$ for any $\eta>0$. So many results on the geometric
and statistical properties of unimodal maps with non--flat
critical point extend to maps which are only $C^{2+\eta}$.

Since we cannot use the negative Schwarzian property at all here,
we must look rather closely at the behaviour of the map on small
scales. We use a result in \cite{ms} to estimate the cross ratio
distortion in terms of sums of lengths of intervals.  We split up
this sum into blocks using the domains of first return maps to
small intervals around the critical point.  The precise behaviour
of the branch containing the critical point, the central branch,
determines how we choose our blocks. Since we have no negative
Schwarzian property, there are particular difficulties when a
block of our sum contains points which spend a very long time in
the central branch (when there is a so called `saddle node
cascade' or an `Ulam-Neumann cascade'). The main tool we use here
is the \emph{real bounds} proved by \cite{v,sbounds,svarg}.
Roughly speaking, these results give us a sequence of first return
maps where the diffeomorphic branches have a uniformly large
extension. This gives bounded distortion of the derivative on
these branches which allows us to estimate the sums of lengths of
intervals.

\subsection{Statement of the main result}
We explain the terminology in the following definitions.  Given an
interval $T$, and a subinterval $J \subset T$, we defined the
cross--ratio $B(T,J)$ above. Note that if we again denote the
left--hand and right--hand components of $T \setminus J$ by $L$
and $R$ respectively, we have another measure of cross--ratio
$$A(T, J) := \frac{|T||J|}{|L\cup J||J \cup R|},$$ (however, we
focus mainly on $B(T, J)$).

Suppose that $g : T \to \RRR$ is a diffeomorphism.  We define
$B(g, T, J)$ as above, but we also have
$$A(g, T, J) := \frac{A(g(T), g(J))}{A(T, J)},$$ another estimate
of how the map distorts cross--ratios.  Observe that for
diffeomorphisms $g:T \to g(T)$ and $h: g(T) \to h\circ g(T)$ we
have
$$B(h\circ g, T, J) = B(h, g(T), g(J)) B(g, T, J).$$  Similarly for
$A(g, T, J)$.

We say that $T$ is a {\em $\delta$--scaled neighbourhood of $J$}
if $\frac{|L|}{|J|}, \frac{|R|}{|J|} > \delta$.  We suppose
throughout that our functions map from $I:=[0,1]$ into itself, and
$\bd I $ into $\bd I$.

We say that a unimodal $C^k$ map $g$ has {\em non--flat} critical
point $c$ if there exists some neighbourhood $U$ of $c$ and a
$C^k$ diffeomorphism $\phi: U \to I$ with $\phi(c) = 0$ \st $g(x)
= \pm|\phi(x)|^\alpha +  g(c)$ for some $\alpha>1$. The value
$\alpha$ is known as the {\em critical order} for $g$. We denote
the set of such maps by $NF^k$ and this neighbourhood by $U_\phi$.

Such maps have many good properties.  For example, they have no
wandering intervals, see for example Chapter IV of \cite{ms}. More
importantly for us here is how such maps distort cross--ratios. In
particular, how iterates of such maps distort cross--ratios.  Our
main result is as follows.

\begin{thm}
For any $\eta>0$, let $f \in NF^{2+\eta}$ be a unimodal map with a
critical point whose iterates do not converge to a periodic
attractor. Then for any $0<K<1$, there is an interval $V$ around
the critical point \st if, for an interval $T$ and some $n>0$,
\begin{itemize}
\item $f^n|_T$ is monotone; and
\item $f^n(T) \subset V$,
\end{itemize}
then $$B(f^n, T, J) > K,$$
$$A(f^n, T, J) > K$$ where $J$ is any subinterval of $T$.
\label{thm:cr}
\end{thm}

This theorem is proved for $C^3$ maps in \cite{kozsch}. Note that
in fact we prove that if $0<\eta\le 1$ then for any
$0<\eta'<\eta$, there exists $C>0$ \st if $J,T,V$ are as in the
theorem then $A(f^n, T, J), B(f^n,T,J)> \exp\{-C(\sup_j
|V_j|)^{\eta'}\}$.

\subsection{Strategy of the proof}
Our setup will involve first return maps to a neighbourhood of
$c$, as outlined below. For the case where $c$ is non--recurrent
see \cite{strien}. So we suppose throughout that $c$ is recurrent.

An open interval $V$ is {\em nice for $f$} if $f^n(\bd V) \cap V =
\OOO$ for $n \ge 1$.  (When it is clear what $f$ is, we just refer
to such interval as nice.) It is easy to see that we can find
arbitrarily small nice intervals around $c$.

Let $I_0\ni c$ be a nice interval.  For every $x\in I$ whose orbit
intersects $I_0$, let $n(x):=\min\{k>0:f^k(x)\in I_0\}$.  If
additionally $x\in I_0$, let $I_0^j\ni x$ be the maximal
neighbourhood \st $f^{n(x)}(I_0^j) \subset I_0$.  We obtain the
first return map $F_0:\bigcup_jI_0^j \to I_0$.  We label the
interval which contains $c$ by $I_0^0$; this interval is called
the \emph{central domain}. Observe that $F_0$ is a diffeomorphism
on all domains $I_0^j$ except when $j=0$.  $F_0$ is unimodal on
$I_0^0$.  Note also that $I_0^0$ is again a nice interval.  We
will call it $I_1$ for the next step in the inducing process; i.e.
we define $F_1:\bigcup_j I_1^j \to I_1$ to be the first return map
to $I_1=I_0^0$.  It has central domain $I_1^0=I_2$.  Continuing
inductively, we obtain maps $F_i:\bigcup_jI_i^j\to I_i$.  The
sequence $I_0 \supset I_1\supset\cdots$ is called the
\emph{principal nest}, and $F_i|_{I_i^j}: I_i^j\to I_i$ is a
\emph{branch} of $F_i$.

If $x \notin I_i$ but $n(x)$ is defined then there is a maximal
interval $U_i^j\ni x$ \st $f^{n(x)}: U_i^j\to I_i$ is a
diffeomorphism. So we may extend $F_i$, letting
$F_i|_{U_i^j}:U_i^j \to I_i$. Then letting $\bigcup_j U_i^j$
consist of all such intervals added to $\bigcup_jI_i^j$, we call
$F_i:\bigcup_j U_i^j \to I_i$ the {\em first entry map to $I_i$}.
We will often switch between these two very similar types of map.

\begin{figure}[htp]
\begin{center}
\includegraphics[width=0.35\textwidth]{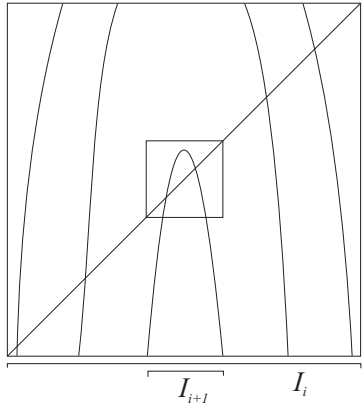}
\caption{$F_i$ is high and central.} \label{fig:central}
\end{center}
\end{figure}

For simplicity, except in the appendix, we will assume that
$F_i(c)$ is a maximum for $F_i|_{I_{i+1}}$. We say that $F_i$ is
{\em low} if $F_i(c)$ lies to the left of $c$ and $F_i$ is {\em
high} if $F_i(c)$ lies to the right of $c$. $F_i$ is {\em central}
if $F_i(c)$ is inside $I_{i+1}$ (if this is not the case, then
$F_i$ is {\em non--central}). Figure~\ref{fig:central} shows $F_i$
which is high and central return.

Suppose that $f^n:T \to f^n(T)$ is a diffeomorphism and $f^n(T)
\subset I_{0}$. It can be shown (see Lemma~\ref{thm:mu}) that we
get a lower bound on $B(f^n, T, J)$ if we can find some bound on
$\sum_{k=0}^{n-1}|f^k(T)|$. In fact, we consider
$\sum_{k=0}^{n-1}|f^k(T)|^{1+\xi}$ for some $0<\xi<\eta$. We will
split up this sum into blocks determined by the principal nest
explained above. Note that our proofs extend easily to $A(f^n, T,
J)$, see \cite{strien}.

We fix $n$ and $T$ as in Theorem~\ref{thm:cr}, let $n_0=n$.  For
$i> 0$, suppose that some iterate $f^j(T)$ enters $I_i$ for $0\le
j\le n$. Now we let $n_i$ be the last time that $f^j(T)\subset
I_i$, i.e. $f^{n_i}(T) \subset I_i$ and $f^{n_i+j}(T) \nsubseteq
I_i, \ 0<j \le n-n_{i}$. If $f^j(T)$ is never contained in $I_i$
for $0\le j\le n$ then we let $n_i = n_{i-1}$. For each $i$, we
will be interested in estimating
$$\sum_{k=1}^{n_i-n_{i+1}} |f^{k+n_{i+1}}(T)|^{1+\xi}\ \ \
\hbox{ we call this the {\em the sum for }} F_i.$$ As we will see
later, if $F_i$ is non--central infinitely often then
Theorem~\ref{thm:svarg} implies that as $i \to \infty$ the
intervals $I_i$ shrink down to $c$. Thus we are able to bound
$\sum_{k=0}^{n-1}|f^k(T)|^{1+\xi}$ by bounding the sums for all
$F_i$.  We will use a slightly different method when there exists
a nice $I_0$ \st $F_i$ is always central.

In order to prove the main theorem, we will consider the following
cases.  Note that we only assume that $f \in NF^2$ in the
following three propositions.
\begin{itemize}
\item $F_{i-2}$ is non--central.  We consider the sum for $F_i$
whenever $f^j(T)\cap \bd I_{i+1} =\es$ for all $0\le j<n_i$, as
follows.
\begin{prop} Suppose that $F_{i-2}$ is non--central and $f^j(T)
\cap \bd I_{i+1} =\es$ for all $0\le j<n_i$. Then there exists
$C_{wb}>0$ \st $$\sum_{k=1}^{n_i-n_{i+1}} |f^{k+n_{i+1}}(T)| <
C_{wb} \sigma_{i} \frac{|f^{n_i}(T)|}{|I_{i}|},$$ where $\sigma_i
:= \sup_{V \in \{I_i^j\}_j }\sum_{k=1}^{n(V)} |f^k(V)|$ (and
$n(V)$ is defined as $k$ where $F_i|_V = f^k$).
\label{prop:entiresum}
\end{prop}
We call this a {\em well bounded case}.  It is dealt with in
Section~\ref{sec:bdd}.

\item $F_{i-2}$ is non--central and $F_i, \ldots, F_{i+m-1}$
are central.  We consider the sums for $F_i, F_{i+1}, \ldots,
F_{i+m}$ whenever $f^j(T)\cap \bd I_{i+m+1} =\es$ for all $0\le
j<n_i$, as follows.
\begin{prop}
Suppose that $F_{i-2}$ is non--central, $F_i, \ldots, F_{i+m-1}$
are central and $f^j(T)\cap \bd I_{i+m+1} =\es$ for all $0\le
j<n_i$.  For all $\xi>0$ there exists $C_{casc}>0$ \st
$$\quad\quad\quad \sum_{k=1}^{n_i- n_{i+m+1}}|f^{ k+n_{i+m+1}}(T)|^{1+\xi}
 < C_{casc}\sigma_{i, m}\max_{n_{i+m+1}< k\le n_i}
|f^k(T)|^\xi$$ \label{prop:entiresum'} where $\sigma_{i, m}$ is
defined as follows.  Let $\sigma_i := \sup_{V \in \{I_i^j\}_j}
\sum_{k=1}^{n(V)} |f^k(V)|$.  Let $\hat{V} \subset I_i \setminus
I_{i+1} $ be an interval \st $f^{\hat{n}}(\hat{V})$ is one of the
connected components of $I_i \setminus I_{i+1} $ for some
$\hat{n}>0$ and $f^j(\hat{V})$ is disjoint from both $I_i
\setminus I_{i+1} $ and $I_m$ for $0< j< \hat{n}(\hat{V})$.  Then
$\sigma_{i, m}$ is the supremum of all such sums
$\sum_{j=1}^{\hat{n}(\hat{V})}|f^j(\hat{V})|$ and $\sigma_i$.
\end{prop}
We call this the {\em cascade case}.  It is dealt with in
Section~\ref{sec:cascade}.

\item $F_{i-2}$ is central and $F_{i-1}$ is high and non--central.  We
consider the sum for $F_i$ whenever $f^j(T)\cap \bd I_{i+1} =\es$
for all $0\le j<n_i$, as follows.
\begin{prop}
Suppose that $F_{i-2}$ is central, $F_{i-1}$ is high and
non--central and $f^j(T)\cap \bd I_{i+m+1} =\es$ for all $0\le
j<n_i$. Then there exist $C_{ex}>0$ and $n_{i+1}< n_{i,3}< n_{i,2}
< n_i$ \st $f^{n_{i,2}}(T), f^{n_{i,3}}(T) \subset I_{i}$ and
$$\quad\quad \sum_{k=1}^{n_i-n_{i+1}} |f^{k+n_{i+1}}(T)| < C_{ex} \sigma_{i}
\left(\frac{|f^{n_i}(T)|}{|I_{i}|} +
\frac{|f^{n_{i,2}}(T)|}{|I_{i}|}+ \frac{|f^{n_{i,3}}(T)|}{|I_{i}|}
\right).$$ \label{prop:entiresumex}
\end{prop}
\noindent (In some cases, the last two terms in the sum are not
required.) We call this the {\em exceptional branches case}.   It
is dealt with in Section~\ref{sec:exceptional}.  We also note
there that if $F_{i-2}$ is central and $F_{i-1}$ is low and
non--central then we are in another well bounded case, and so the
conclusion of Proposition~\ref{prop:entiresum} holds.

\item We have an interval $I_0$ \st $F_i$ are all central for
$i=0, 1, \ldots$.  We call this the {\em infinite cascade case}.
We prove Theorem~\ref{thm:cr} for this case in
Section~\ref{sec:infcasc}.
\end{itemize}

The proof of Theorem~\ref{thm:cr} for the non--infinite cascade
case is given in Section~\ref{sec:proof of cr non-inf}.

With these propositions, for $0<\eta'<\eta$, we can decompose the
sum $\sum_{k=0}^{n-1}|f^k(T)|^{1+\eta'}$ into blocks of sums
$\sum_{k=1}^{n_i-n_{i+1}} |f^{k+n_{i+1}}(T)|^{1+\eta'}$.  We then
show that each of these is uniformly bounded.  We will then show
that $\sum_{k=1}^{n_i-n_{i+1}} |f^{k+n_{i+1}}(T)|^{1+\eta}$ decays
in a uniform way with $i$.

The first two cases use real bounds of Theorem~\ref{thm:svarg}.
These bounds imply that $B(I_i^j,I_i)$ are bounded above. This
will also be true for all except possibly two domains of $F_i$ in
the third case.  The main tool here is Lemma~\ref{lem:lambda},
which gives us some decay of cross--ratios when we have these real
bounds. Note that the conditions $f^j(T)\cap \bd I_{i+1} =\es$ for
all $0\le j<n_i$ in well bounded and exceptional cases, and
$f^j(T)\cap \bd I_{i+m+1} =\es$ for all $0\le j<n_i$ in the
cascade case, make the propositions simpler to prove. However, as
we remark in Section~\ref{sec:proof of cr non-inf}, it is easy to
see how to split up the intervals in the other cases in order to
prove Theorem~\ref{thm:cr}.

The final case, which arises in the infinitely renormalisable
case, is different from the other three.  We use a lemma of
\cite{kozsch} to find some uniform expanding property which helps
bound the sums.

In all cases except the infinite cascade case we must ensure that
we have some initial interval which has a first return map which
is well bounded. To do this we can simply pick some nice interval
to begin with and then induce until we find a map which is
well--bounded. This is always possible when there is not an
infinite cascade.

Note that we need extra smoothness to bound cross--ratios in the
cascade case. This ensures that we can deal with the case when we
have many consecutive low central returns, a `saddle node
cascade'.

In Kozlovski's proof for $C^3$ maps he was able to use the fact
that there exists some $C>0$ depending only on $f$ \st for
interval $J \subset T$ we have $B(f, T, J)
> \exp\{-C|T|^2\}$ and $A(f, T, J) >
\exp\{C|L||R|\}$. See Chapter IV.2 of \cite{ms}.  In particular
this means that there exist such real bounds as in
Theorem~\ref{thm:svarg} for all $i$, not just those for which
$F_{i-1}$ is a non--central return. So the long central cascades
we encounter in Section~\ref{sec:cascade} present much less of a
problem in the $C^3$ case.  Indeed, the work done in
Section~\ref{sec:exceptional} is also unnecessary in the $C^3$
case.

We will deal with the well bounded case first. It is the simplest
and gives us a good idea about how we may proceed in general.  We
will use $J$ to refer to a general interval from here until
Section~\ref{sec:proof of cr non-inf}.  This allows us to use less
notation. When we use the constant $C>0$, we mean some constant
depending only on $f$.

%Our results answer a question of Shen, posed at the December
%meeting of the University of Warwick Dynamics Symposium 2002-2003
%which arose out of his work as well as that of van Strien and
%Vargas.

{\bf Acknowledgements}

This work was undertaken as part of my thesis at the University of
Warwick, which was funded by the EPSRC. I would like to thank my
supervisor Oleg Kozlovski for his support. I would also like to
thank the dynamical systems group at Warwick, in particular
Weixiao Shen whose suggestions, support and enthusiasm were
invaluable.  Further thanks to Sebastian van Strien, Henk Bruin
and to the referee for useful comments.

\section{Introductory results}
\label{sec:basiccrstuff} Without loss of generality, we suppose
throughout that our maps have a maximum at the critical point. We
also suppose that $f$ is symmetric about $c$. That is,
$f(c-\epsilon)=f(c+\epsilon)$ for all $\epsilon$. This assumption
is useful for simplifying proofs (particularly in
Section~\ref{sec:exceptional}, which is already quite technical),
but is not crucial since on small scales our maps will be
essentially symmetric (in particular, $|Df(c-\epsilon)|$ and
$|Df(c+\epsilon)|$ are arbitrarily close for small enough
$\epsilon$). We let $C\le |g|_U\le C'$ mean $\sup_{x\in
U}|g(x)|\le C'$ and $\inf_{x\in U}|g(x)|\ge C$.

The following theorem is proved for a more general case in Chapter
IV of \cite{ms}.  Here we will let $w_g$ be the modulus of
continuity of a continuous map $g$, i.e.
$w_g(\epsilon):=\sup_{|x-y|<\epsilon}|g(x)-g(y)|$.

\begin{thm}
For a unimodal map $g:I \to I$, $g \in NF^2$, if $T$ is an
interval \st $g^n|_T$ is a diffeomorphism and $J \subset T$ is a
subinterval, then there exists some $C>0$  \st
$$B(g^n, T, J) > \exp\left\{-C  \sum_{i=0}^{n-1}w_{D^2g}(|g^i(T)|)
|g^i(T)| \right\}.$$ This bound also holds for $A(f^n, T, J)$.
\label{thm:mu}
\end{thm}

In Sections~\ref{sec:proof of cr non-inf} and \ref{sec:infcasc} we
will use the fact that when $g\in NF^{2+\eta}$ for some $\eta>0$,
we can replace $Cw_{D^2g}(\epsilon)$ by $C\epsilon^\eta$.

The following lemma, a consequence of the absence of wandering
intervals, is Lemma 5.2 in \cite{kozsch}.

\begin{lem}
Suppose that $g \in NF^2$, $g:I \to I$. Then there exists a
function $\tau:[0, |I|] \to [0, \infty)$ \st $\lim_{\epsilon \to
0}\tau(\epsilon) = 0$ and for any interval $V$ for which $g^n|_V$
is a diffeomorphism and $g^n(V)$ is disjoint from the immediate
basins of periodic attractors, we have
$$\max_{0 \le i \le n}|g^i(V)| < \tau(|g^n(V)|).$$
\label{lem:tau}
\end{lem}

We may use this lemma and Theorem~\ref{thm:mu} to get
\begin{equation} B(g^n, T, J) > \exp\left\{-\sigma'(|g^{n-1}(T)|)
\sum_{i=0}^{n-1} |g^i(T)| \right\} \label{eqn:sigma'sum}
\end{equation} whenever $f^n(T)$ is disjoint from the immediate
basins of periodic attractors, where
\begin{equation} \sigma'(|g^m(T)|) = Cw_g\circ
\tau(|g^m(T)|).\label{eqn:sigma'}\end{equation}

We will use the following result of \cite{svarg} throughout.  (In
fact it is stated there in greater generality, as Theorem A.)

\begin{thm}
If $g \in NF^2$ is a unimodal map with recurrent critical point,
then the following hold.
\begin{itemize}
\item[(a)] For all $k \ge 0$ there exists $\xi(k)
> 0$ \st if $G_{i-1}: \bigcup_j I_{i-1}^j \to I_{i-1}$ is
non-central, then $I_{i+k}$ is a $\xi(k)$--scaled neighbourhood of
$I_{i+k+1}$.

\item[(b)] For each $\xi>0$ there is some $\hat\xi>0$ \st if $I_i$ is a
$\xi$--scaled neighbourhood of $I_{i+1}$ then $I_{i+1}$ is a
$\hat\xi$--scaled neighbourhood of any domain of $G_{i+1}$.
\end{itemize}
\label{thm:svarg}
\end{thm}

This result gives us {\em real bounds} for some of our first
return maps. We let $\chi:= \xi(1)>0$ from the above theorem for
our map $f$.

The following theorem is an improvement of the classical Koebe
Lemma. It is presented in more generality in \cite{svarg} as
Proposition 2: `a Koebe principle requiring less disjointness'.
Note that actually for our purposes, the classical Koebe Lemma is
enough.

\begin{thm}
Suppose that $g \in NF^2$.  Then there exists a function $\nu: [0,
|I|] \to [0, \infty)$ \st $\nu(\epsilon) \to 0$ as $\epsilon \to
0$ with the following properties.  Suppose that for some intervals
$J \subset T $ and a positive integer $n$ we know that $g^n|_T$ is
a diffeomorphism.  Suppose further that $g^n(T)$ is a
$\delta$--scaled neighbourhood of $g^n(J)$ for some $\delta>0$.
Then,
\begin{itemize}
\item[(a)] for every $x, y \in J$,
$$\frac{|Dg^n(x)|}{|Dg^n(y)|} < \exp\left\{\nu(S(n, T))
\sum_{i=0}^{n-1}|g^i(J)|\right\}
\left[\frac{1+\delta}{\delta}\right]^2=:C(\delta)$$ where
$S(n,T):=\max_{0\le k\le n-1}|f^k(T)|$.

\item[(b)]  $T$ is a
$\tilde\delta$--scaled neighbourhood of $J$ whenever
$$\tilde\delta := \frac{1}{2}\exp\left\{-\theta\right\}
\left[\frac{1+\delta}{\delta}\right]^2 \left(\frac{-2\theta
+\delta(1-2\theta)}{2+\delta}\right)$$ is positive, where
$\theta:=\nu(S(n, T)) \sum_{i=0}^{n-1}|g^i(J)|$.
\end{itemize}
\label{thm:disjkobcr}
\end{thm}

Again we may use Lemma~\ref{lem:tau} to substitute $\nu(S(n, T))$
with $\nu'(|f^n(T)|)$ where we define $\nu'(|f^m(V)|) :=
\nu\circ\tau(|f^m(V)|)$.  We will use the result of
Theorem~\ref{thm:svarg}(b) extensively, but we use $\tilde\delta$
when $\theta=\nu'(|I_0|)$. Usually $\delta$ will be related to the
$\chi$ we obtained following Theorem~\ref{thm:svarg}.

We will sometimes be in a situation where we wish to estimate the
derivative of a function in between two points at which we know
something about the derivative. The following two well known
results allow us to do this. The following is known as the {\em
Minimum Principle}; see, for example, Theorem IV.1.1 of \cite{ms}.

\begin{thm}
Let $T=[a, b] \subset I$ and $g:T \to g(T) \subset I$ be a $C^1$
diffeomorphism.  Let $x \in (a, b)$.  If for any $J^* \subset T^*
\subset T$,
$$B(g, T^*, J^*) > \mu_g >0$$
then $$|Dg(x)| > \mu_g^3 \min(|Dg(a)|, |Dg(b)|).$$ \label{thm:min}
\end{thm}

To see a proof of the following well known result see \cite{ms}.

\begin{thm}
For $g \in NF^2$ there exist $n_0 \in \NNN$ and $\rho_g
>1$ \st if $p$ is a periodic point of period $n \ge n_0$ then
$|Dg^n(p)|>\rho_g.$ \label{thm:B}
\end{thm}

We are now ready to begin the proof of Theorem~\ref{thm:cr}.

\section{Well bounded case}
\label{sec:bdd}

Here we deal with the case where $F_{i-2}$ is non--central and
$f^j(T)\cap \bd I_{i+1} =\es$ for all $0\le j<n_i$.  In our
estimates, we are principally interested in iterates of $T$
landing in $I_i^j$ for $j\neq 0$.  By Theorem~\ref{thm:svarg}, the
fact that $F_{i-2}$ is non--central implies that the first return
domains $I_i^j$ are all well inside $I_i$.  This enables us to
estimate the sum for $F_i$, and is the reason we call this case
well bounded.

Let $n_i' > n_{i+1}$ be minimal \st $f^{n_i'}(T) \subset I_i$.  We
will initially assume that we have some $\kappa>0$ \st for the
`return sum',
\begin{eqnarray} \sum_{k=0}^{j_i}|F_i^{k}(f^{n_i'}(T))| < \kappa
|f^{n_i}(T)| \label{eqn:kappa} \end{eqnarray} where $j_i$ is \st
$F^{j_i}|_{f^{n_i'}(T)} = f^{n_i - n_i'}|_{f^{n_i'}(T)}$. We prove
Proposition~\ref{prop:entiresum} before bounding this return sum
in order to give an idea why we need bounds on return sums. Except
for the proof of \eqref{eqn:kappa}, this is similar to the proof
of Lemma 5.3.4 of \cite{kozthe}. There, it is assumed that $f \in
C^3$ in order to bound the sum
$\sum_{k=0}^{j_i-1}|F_i^{k}(f^{n_i'}(T))|$. Those methods fail in
the $C^2$ case.

\begin{proof}[ Proof of Proposition~\ref{prop:entiresum} assuming
(\ref{eqn:kappa})] Let $n_{i+1}=m_0 < m_{1} < \cdots < m_{j_i} =
n_{i}$ be all the integers between $n_{i+1}$ and $n_{i}$ \st
$f^{m_{j}}(T) \subset I_i \setminus I_{i+1}$ for $j=1, \ldots,
j_i-1$ and let $m_{0} = n_{i+1}$. Now let $F_i:\bigcup_j U_i^j \to
I_i$ be the first entry map to $I_i$. We will decompose
$\sum_{k=1}^{n_i-n_{i+1}} |f^{k+n_{i+1}}(T)|$ as
$\sum_{j=0}^{j_i-1} \sum_{k=1}^{m_{j+1} - m_{j}}
|f^{k+m_{j}}(T)|$.

For $1 \le j \le j_i -1$ and $1\le k<m_{j+1}-m_j$, let $U_i^l$ be
the domain of first entry to $I_i$ \st $f^{m_{j}+k}(T) \subset
U_i^l$. Suppose that $F_i|_{U_i^l}=f^{i_l}$. Then there exists an
extension to $V_i^l \supset U_i^l$ so that $f^{i_l}: V_i^l \to
I_{i-1}$ is a diffeomorphism. Then by the Koebe Lemma we have the
distortion bound: $\frac{|f^{k+m_j}(T))|}{|U_i^l|} \le C(\chi)
\frac{|f^{m_{j+1}}(T)|}{|I_i|}$.  Whence
\begin{align*}
\sum_{k=1}^{m_{j+1} - m_{j}} |f^{m_{j}+ k}(T)| & \le C(\chi)
\left(\frac{|f^{m_{ j+1}}(T)|}{|I_{i}|}\right) \sum_{k=0}^{m_{j+1}
- m_j -1} |f^k(U_i^j)| \\
& \le C(\chi) \sigma_{i}
\frac{|f^{m_{j+1}}(T)|}{|I_{i}|}.\end{align*} Therefore
$$\sum_{k=1}^{n_i-n_{i+1}} |f^{k+n_{i+1}}(T)| \le C(\chi)
\frac{\sigma_{i}}{|I_{i}|} \sum_{j=1}^{j_i}|f^{m_{j}}(T)|= C(\chi)
\frac{\sigma_{i}}{|I_{i}|} \sum_{k=0}^{j_i-1}|F_i^k(\hat{T})|$$
where $\hat{T} := f^{n_{i}'}(T)$. This is bounded above by
$\kappa|f^{n_i}(T)|$ due to (\ref{eqn:kappa}), so we are finished.
\end{proof}

\subsection{Bounding return sums}
\label{subsec:returns}

In this subsection we will introduce some tools which we use
extensively in the remainder of this paper. We then use these
tools to prove that (\ref{eqn:kappa}) holds.

The proof of the following simple lemma is left to the reader.

\begin{lem}
For all $\delta>0$ there exists $\Delta= \Delta(\delta) >0$ \st
$\Delta(\delta) \to 0$ as $\delta \to \infty$ with the following
property. Suppose that $U$ is an interval, $J \subset U$ is a
subinterval and that the left and right components of $U \setminus
J$ are denoted by $L$ and $R$ respectively. Suppose further that
$|L|, |R| > \delta|J|$. Then
$$B(U, J) < \Delta.$$
\label{lem:Delta}
\end{lem}

Let $D_1$ denote the set of non--central domains  $F_i^{-1}(I_i)$,
i.e. $D_1 = \bigcup_{j \neq 0} I_i^j$.  Let $D_2$ denote the set
of domains $F_i^{-1}(D_1)$ which are disjoint from the central
domain. Inductively, we let $D_k$ denote the set of domains
$F_i^{-1}(D_{k-1})$ which are disjoint from the central domain.
Then for any element $J_k \in D_k$, $F_i^k: J_k \to I_i$ is a
diffeomorphism. We will bound $\sum_{j=0}^{k-1}|F_i^j(J_k)|$ for
any $J_k \in D_k$ by showing that there exists some $\lambda < 1$
independent of $i$ \st for $k>1$ we have $B(I_i , J_k) \le \lambda
B(I_i, F_i(J_k))$. We let
\begin{equation} \mu := \exp\left\{-\sigma'(|I_0|)\right\}
\label{eq:mu}\end{equation} where $\sigma'$ is given by
\eqref{eqn:sigma'}.  By \eqref{eqn:sigma'sum}, if $J', f(J'),
\ldots, f^m(J')$ is a disjoint set of intervals and $J \supset
J'$, we have $B(f^m, J', J)> \mu$. Therefore, if $n(j)$ is the
return time of $I_i^j$ to $I_i$ and $J \subset I_i^j$ then
$B(f^{n(j)}, I_i^j, J)> \mu$.

The following lemma is Lemma 2.3 of \cite{gk}.

\begin{lem}
For every $\delta>0$ there exists $\lambda' = \lambda'(\delta)<1$
\st if $J \subset V \subset U$ are intervals and $U$ is a
$\delta$--scaled neighbourhood of $V$ then
$$B(U, J) < \lambda' B(V, J).$$
Furthermore, $\lambda' \to 1$ as $\delta \to 0$.
\label{lem:lambda'}
\end{lem}

We add this lemma to \eqref{eqn:sigma'sum} as follows.

\begin{lem}
Given $\delta>0$, there exist  $0< \lambda=\lambda(\delta)<1$ and
$\epsilon>0$ \st if $|I_0|<\epsilon$ and $I_{i-1}$ is a
$\delta$--scaled neighbourhood of $I_i$, then for any $J \subset
I_i^j$ with $j \neq 0$,
$$B(I_i, J) < \lambda B(I_i, F_i(J)).$$
\label{lem:lambda}
\end{lem}
\begin{proof}
From the previous lemma there exists some $\lambda' =
\lambda'(\delta) <1$ \st
$$B(I_i, J) < \lambda' B(I_i^j, J).$$   Now from
\eqref{eqn:sigma'sum} we obtain
$$B(I_i, J) < \lambda' \frac{ B(I_i, F_i(J))}{\mu}$$ where $\mu$
is defined in \eqref{eq:mu}.  Since $\mu \to 1$ as $|I_0| \to 0$,
if $\epsilon$ is chosen small enough then $\frac{\lambda'}{ \mu} <
1$. We let $\lambda := \frac{\lambda'}{ \mu}$. Thus $B(I_i, J) <
\lambda B(I_i, F_i(J))$.
\end{proof}

We will consider $\lambda = \lambda(\tilde\chi)$ where
$\tilde\chi$ comes from Theorem~\ref{thm:disjkobcr}(b) applied to
$\chi$ and $\chi$ comes from Theorem~\ref{thm:svarg}(a), i.e.
$\tilde\chi$ takes the role of $\delta$ in Lemma~\ref{lem:lambda}.
In fact we shall adjust $\lambda$ again in
Section~\ref{sec:exceptional}, but it will remain independent of
$i$ and strictly less than 1.

\begin{proof}[Proof of \eqref{eqn:kappa}.]
For $k \ge 2$, $B(I_i, J_k) < \lambda^{k-1} B(I_i,
F_i^{k-1}(J_k))$.  Suppose that $F_i^{k-1}(J_k) \subset I_i^j$.
Then by Lemma~\ref{lem:Delta}, using Theorems~\ref{thm:svarg} and
\ref{thm:disjkobcr} (b), $B(I_i, I_i^j)< \Delta$ where $\Delta=
\Delta(\tilde\chi)$.  Thus, it is easy to see $B(I_i,
F_i^{k-1}(J_k)) < \Delta \frac{|F_i^{k-1}(J_k)|}{|I_i^j|}$.  Now
by the Koebe Lemma, $|F_i^{k-1}(J_k)| < C(\chi) |F_i^{k}(J_k)|
\frac{|I_i^j|}{|I_i|}$, so we know that $B(I_i, F_i^{k-1}(J_k)) <
C(\chi) \Delta \frac{|F_i^{k}(J_k)|}{|I_i|}$. We apply these
estimates to the sizes of $J_k$: $$|J_k|< \frac{
|I_i|}{1+\frac{2|I_i|}{\lambda^{k-1} C(\chi) \Delta |F_i^{k}(J_k)|
}}.$$ Then $ |J_k| < C \lambda^{k-1}|F_i^{k}(J_k)|$. So
$\sum_{j=0}^{k-1} |F_i^j(J_k)| < C \frac{|F_i^{k}(J_k)|}{1-
\lambda}$. Whence $$\sum_{j=0}^{k} |F_i^j(J_k)| <
|F_i^{k}(J_k)|\left(1+ \frac{C}{1-\lambda}\right).$$ This holds
for any sum of returns which never lands in the central domain. It
is independent of $i$. Letting $\kappa=\left(1+
\frac{C}{1-\lambda}\right)$ we prove (\ref{eqn:kappa}).
\end{proof}

\section{Cascade case}
\label{sec:cascade}

This section is devoted to the proof of
Proposition~\ref{prop:entiresum'}.  Note that if there is a
uniform upper bound on the length of sequences  $F_i, F_{i+1},
\ldots, F_{i+m}$ all having central returns then
Theorem~\ref{thm:svarg} implies that we may prove
Proposition~\ref{prop:entiresum'} as a well bounded case. However,
there may be arbitrarily long sequences of consecutive central
returns.

\begin{proof}[Proof of Proposition~\ref{prop:entiresum'}]
We suppose that there $i$ is \st $f^{n_i}(T) \subset I_i$ where
$F_{i-2}$ has a non--central return and $F_{i+j}$ all have central
returns for $j=0,\ldots, m-1$ and that $F_{i+m}$ has a non-central
return. For $\xi>0$ we will bound the sum
$$\sum_{k=1}^{n_i- n_{i+m+1}}|f^{ k+n_{i+m+1}}(T)|^{1+\xi}.$$
For our intial estimates, we may omit $\xi$, but later it will be
necessary to include it.  Recall that we always assume here that
$f^j(T)\cap \bd I_{i+m+1} =\es$ for all $0\le j<n_i$.

Let $m_0= n_{i+m+1}$ and let  $m_0<m_1\le n_i$ be the smallest
integer \st $f^{m_1}(T) \subset I_i \setminus I_{i+1}$. Let
$m_1<m_2\le n_i$ be the next integer for which $f^{m_2}(T) \subset
I_i \setminus I_{i+1}$ if such $m_2$ exists. Proceeding in this
manner, we obtain a sequence, $n_{i+m+1} < m_1 < m_2 < \cdots <
m_N = n_i$.

So  $$\sum_{k=1}^{n_i- n_{i+m+1}} |f^{k+n_{i+m+1}}(T)| =
\sum_{j=0}^{N-1} \sum_{k=1}^{m_{j+1}- m_{j}} |f^{k+ m_{j}}(T)|.$$

Define $m_{N-1}<m'\le n_i$ to be minimal \st $f^{m'}(T)\subset
I_i\sm I_{i+m+1}$.  Assuming that $F_i|_{I_i^0} =f^s$, there
exists $0\le p\le m$ such that $m'+sp=m_{N}=n_i$.  We can rewrite
the sum
\begin{align*}
\sum_{k=1}^{n_i- n_{i+m+1}} |f^{k+n_{i+m+1}}(T)| & =
\sum_{j=0}^{N-2} \sum_{k=1}^{m_{j+1}- m_{j}} |f^{k+
m_{j}}(T)|+\sum_{k=1}^{m'-m_{N-1}}|f^{k+m_{N-1}}(T)|\\
&\quad + \sum_{r=0}^{p-1}\sum_{k=1}^s|f^{k+rs+m'}(T)|.
\end{align*}

Using the method from the proof of
Proposition~\ref{prop:entiresum},
$$\sum_{k=1}^{m'}|f^{k+m_{N-1}}(T)|+
\sum_{r=0}^{p-1}\sum_{k=1}^s|f^{k+rs+m'}(T)|\le
C(\chi)\frac{\sigma_{i,m}}{|I_i|}\sum_{r=0}^{p}|f^{rs+m'}(T)|.$$
We will deal with the sum on the right hand side later.  We will
first show that $\sum_{j=0}^{N-2} \sum_{k=1}^{m_{j+1}- m_{j}}
|f^{k+ m_{j}}(T)| \le C\sigma_{i,m}\frac{|f^{m'}(T)|}{|I_i|}$.

We denote the left and right components of $I_j \setminus I_{j+1}$
by $L_j$ and $R_j$ respectively.  We know from
Theorem~\ref{thm:svarg}(a) and (b) that $\frac{|L_i|}{|I_{i+1}|},
\frac{|R_i|}{|I_{i+1}|} > \hat\chi$.

We define $\hat F_i:\bigcup_j\hat I_i^j \to I_i\sm I_{i+1}$ to be
the first return map to $I_i \sm I_{i+1}$, such that $\hat
F_i(\hat I_i^j) \in \{L_i,R_i\}$. As in the well bounded case, for
each $1\le j\le N-2$ and $1\le k<m_{j+1}-m_j$, there exists a
first entry domain $\hat U$ to $I_i\setminus I_{i+1}$ \st
$f^{k+m_j}(T)\subset \hat U$.  We may assume that
$f^{m_{j+1}-m_j-k}(\hat U)=L_i$. Indeed, for $1\le j\le N-3$ there
exists $\hat I_i^l$ \st $f^{m_{j+1}}(T)\subset \hat I_i^{l}\subset
L_i$. We show that $\hat I_i^{l}$ is well inside $L_i$, which will
allow us to estimate $\frac{|f^{k+m_j}(T)|}{|\hat U|}$.

Suppose that $F_i|_{\hat I_i^{l}}=f^{i_{l}}$.  Then there exists
an extension to $V_i^{l}\supset \hat I_i^{l}$ \st
$f^{i_{l}}:V_i^{l} \to I_{i-1}$.  Clearly $V_i^{l} \subset L_i$,
otherwise niceness is contradicted.  By
Theorems~\ref{thm:svarg}(a) and \ref{thm:disjkobcr}(b), $V_i^{l}$
(and thus $L_i$) is a $\tilde\chi$--scaled neighbourhood of $\hat
I_i^{l}$.

For $1\le j\le N-2$, we have $B(L_i,f^{m_{j+1}}(T))>\mu B(\hat
U,f^{k+m_j}(T))$ where $\mu$ is defined in \eqref{eq:mu}.
Therefore,
$$|f^{k+m_j}(T)|<\frac{|\hat U|}{1+\frac{\mu}{B(L_i,f^{m_{j+1}}(T))}}.$$
As in the well bounded case, using a small adaptation of
Lemma~\ref{lem:lambda}, replacing $F_i$ by $\hat F_i$, we can show
that $B(L_i,f^{m_{j}}(T))<\lambda^{N-1-j} B(L_i,f^{m_{N-1}}(T))$
for $0\le j\le  N-2$.  (Note that $\lambda$ is still the
$\lambda(\tilde\chi)$ discussed following Lemma~\ref{lem:lambda}.)
Therefore, it can be shown that
$$\sum_{j=0}^{N-2} \sum_{k=1}^{m_{j+1}- m_{j}} |f^{k+ m_{j}}(T)|
\le \frac{C\sigma_{i,m}}{1-\lambda}B(L_i,f^{m_{N-1}}(T)).$$ But
since $f^{m_{N-1}}(T)\subset I_i^{j'}$ for some $j'\neq 0$, we
have
$$B(L_i,f^{m_{N-1}}(T))<B(L_i,I_i^{j'})\frac{|f^{m_{N-1}}(T)|}{|I_i^{j'}|}.$$
Notice that $F_i(f^{m_{N-1}}(T))=f^{m'}(T)$.  So the Koebe Lemma
and Lemma~\ref{lem:Delta} give $B(L_i,f^{m_{N-1}}(T))<
C(\chi)\Delta\frac{|f^{m'}(T)|}{|I_i|}$, whence
$$\sum_{j=0}^{N-2} \sum_{k=1}^{m_{j+1}- m_{j}} |f^{k+ m_{j}}(T)|
\le C\sigma_{i,m} \frac{|f^{m'}(T)|}{|I_i|}.$$

It remains to bound $\sum_{r=0}^{p}|f^{rs+m'}(T)|^{1+\xi}$ (as can
be seen below, we only really need $\xi>0$ for our estimates in
the low case). We assume that $f^{m'}(T) \cap \bd I_{i+j} \neq
\OOO$ for $1 \le j <m$: otherwise we have
$\sum_{r=0}^{p}|f^{rs+m'}(T)|^{1+\xi}<|I_i|^{1+\xi}$, and we are
finished.

Let $\hat T = f^{m'}(J)$.  There exists some $M \ge 0$ \st
$F_i^M(\hat{T})= f^{n_i}(T)$.  We will bound
$\sum_{k=0}^M|F_i^k(\hat{T})|^{1+\xi}$.

\begin{figure}[htp]
\begin{center}
\includegraphics[width=0.6\textwidth]{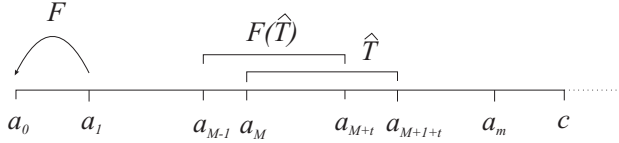}
\caption{When $\hat{T}$ intersects the boundary points $\bd I_j$.}
\label{fig:badcasc}
\end{center}
\end{figure}

If $M$ was uniformly bounded then we would be able to find some
bound on $\sum_{k=0}^M|F_i^k(\hat{T})|$ easily. But $M$ may be
very large.  We consider this sum in two cases: either $F_i$ is
high, or $F_i$ is low (the high case is the most straightforward).
For some background on this dichotomy see \cite{ly}. In both
cases, we relabel $F_i|_{I_{i+1}}$ as $F$ and $I_i$ as $I_0$. Now
let $I_k = (a_k, a_k')$. We are assuming that $F(c)$ is a maximum
for $F$, see Figure~\ref{fig:badcasc}.

{\bf The high case}

We have two cases to consider. We first assume that $F_j$ are high
and central for $j=0, \ldots, m$.  This is known as an
Ulam--Neumann cascade.

\begin{lem}
In the high case, $\sum_{k=0}^M|F^k(\hat{T})|< C|I_0|$.
\label{lem:highcasc}
\end{lem}

\begin{proof}
We know that $I_0$ is a $\hat\chi$--scaled neighbourhood of $I_1$.
We will use the Minimum Principle (Theorem~\ref{thm:min}) and
Theorem~\ref{thm:B} to estimate derivatives.  The idea here is
that either we have derivative uniformly greater than one in
$(a_1, a_m)$ and we can bound $\sum_{k=0}^M|F^k(\hat{T})|$ as a
geometric sum; or we have a small derivative in some region, in
which case we find a bound on the number of $a_i$ that are in this
region.

Let $\gamma>1$ satisfy $\frac\gamma{\gamma-1}>\frac 1{2\hat\chi}$.
Then we may fix some integer $r\ge 1$ \st $2\hat\chi \sum_{i=0}^r
\gamma^{-i}>1$.  Note that $r$ only depends on $\hat\chi$. Observe
that there is a fixed point $p \in (a_1, c)$.  We can choose $I_0$
to be so small that the return time to it is greater than $n_0$
given in Theorem~\ref{thm:B}. Therefore, by that theorem,
$|DF(p)|> \rho_f$. If $|DF(a_1)| \ge \gamma$ then from the Minimum
Principle, $|DF|_{(a_1, p)}> \gamma'$ where $\gamma' =
\mu^3\min(\gamma, \rho_f)$ where $\mu$ is defined in terms of
$|I_0|$ in \eqref{eq:mu}. We fix $I_0$ to be small enough so that
$\gamma'> 1$. Therefore, we have $\sum_{k=0}^M|F^k(\hat{T})|<
\frac {\gamma'}{\gamma'-1} |F_i^M(\hat{T})|$.

Suppose now that there is some $u \in (a_1, c)$ \st $|DF|_{(a_1,
u)}< \gamma$.  We will show that this must mean that $u \in (a_1,
a_r)$ and thus we can uniformly bound the sum of times that $\hat
T$ lies in this region.

Suppose that $(a_1, a_s) \subset (a_1,u)$.  Then we have
$|a_{i+1}- a_i|> \frac{|a_{i}-a_{i-1}|}{\gamma}$ for all $i \le
s-1$. Therefore, if $(a_1, a_s) \subset U$ then
$$|c-a_0| > \sum_{i=0}^{s-1}| a_{i+1}- a_i| > |a_1-a_0| \sum_{i=0}^s
\gamma^{-i}.$$   We know that  $|a_1-a_0| > 2\hat\chi |c-a_0|$. By
the definition of $\gamma$ we must have $s \le r$. Moreover, we
have $|DF|_{(a_s, p)}>  \gamma'$.

This helps us bound $\sum_{k=0}^M|F^k(\hat{T})|$ where
$F^k(\hat{T}) \subset I_0 \setminus I_m$. We suppose that
$F^M(\hat{T}) = (a_0, a_t)$ for $t \le  m$.   See
Figure~\ref{fig:badcasc}.  Then
\begin{eqnarray*} \sum_{k=0}^M|F^k(\hat{T})| & = &
|a_1-a_0|+ \min(2, M-1)|a_2-a_1| + \cdots \\
&  &  + \min(i, M-(i-1))|a_i-a_{t-i}| + \cdots +
|a_{M+t}-a_{M+t-1}|.\end{eqnarray*} This is bounded above by
$$r|a_r-a_0| + |a_N-a_{N+1}| \sum_{i=0}^\infty \frac{\min(i,
M-(i-1)) }{\gamma'^i}.$$ The first summand is bounded by $r|I_0|$
and the second summand is bounded above by $C|a_N-a_{N+1}|$ for
some $C>0$. So we get $\sum_{k=0}^M|F_i^k(\hat{T})|< C|I_0|$ as
required.
\end{proof}

{\bf The low case}

We assume that we are in the same setting as above, but with $F_0$
central and low.  This is known as a saddle node cascade.  Again
we would like to bound $\sum_{k=0}^M|F^k(\hat{T})|$ defined as
above.  However, as we shall see, we are only able to bound
$\sum_{k=0}^M|F^k(\hat{T})|^{1+\xi}$.

\begin{lem}
In the low case, $\sum_{k=0}^M|F^k(\hat{T})|^{1+\xi}<
C|I_0|^{1+\xi}$.\label{lem:lowcasc}
\end{lem}

\begin{proof} We will apply the following result, a form of the Yoccoz
Lemma, see for example \cite{dfdm}.

\begin{lem}
Suppose that $f \in NF^2$.  Then for all $\delta, \delta'>0$ there
exists $C>0$ \st if $I_0$ is a nice interval \st \begin{enumerate}
\item $I_0$ is a $\delta$--scaled neighbourhood of $I_1$; \item
$F_k$ is low and central for $k=0, \dots, m$; \item there is some
$0<k<m$ with $ \frac{|I_k|}{|I_{k+1}|} <1+\delta'$,
\end{enumerate} then for $1 \le k <m$,
$$\frac{1}{C} \frac{1}{\min(k, m-k)^2} < \frac{|I_{k-1}\setminus
I_{k}|}{|I_0|} < \frac{C}{\min(k, m-k)^2}.$$ \label{lem:yoc}
\end{lem}

This lemma was suggested by Weixiao Shen.  For the proof, see the
appendix.  (For comparison with other statements of the Yoccoz
Lemma, note that we will prove that one consequence of our
conditions for the lemma is that we have a lower bound on
$\frac{|I_{m}\setminus I_{m+1}|}{|I_0|}$.)

Suppose that $I_0$ satisfies all the conditions of
Lemma~\ref{lem:yoc}. In particular we assume that for some fixed
$\delta'>0$, we have $\frac{|I_k|}{|I_{k+1}|} <1+\delta'$ for some
$0<k<m$.  Then for any $\xi>0$,

\begin{align*}
\sum_{k=0}^M & |F^k(\hat{T})|^{1+\xi}  \\
& < \sum_{k=0}^m\left( \frac{C|I_0|}{\min(k+t, m-(k+t))^2}+ \cdots
+ \frac{C|I_0|}{\min(k+1,
m-(k+1))^2}\right)^{1+\xi}\\
& < C|I_0|^{1+\xi} \sum_{k=0}^m\left( \frac{1}{k+1}-
\frac{1}{k+t}\right)^{1+\xi}.
\end{align*}
The sum above is bounded above for any $\xi>0$.

Next we suppose that the hypotheses of Lemma~\ref{lem:yoc} do not
hold. In particular, this means $\frac{|I_k|}{|I_{k+1}|}\ge
1+\delta'$ for $k =0, \ldots, m$. Note that $|I_0|\ge
(1+\delta')|I_1| \ge (1+\delta')^2|I_2|\ge \cdots \ge
(1+\delta')^{M}|I_M|$.  Therefore
$$\sum_{k=0}^M|F^k(\hat{T})| <
\frac{1}{2} \sum_{k=0}^M k|I_k|\le \frac{|I_0|}{2} \sum_{k=0}^M
\frac{k}{(1+\delta')^{k}}< C|I_0|.$$ So the lemma is proved.
\end{proof}
We have shown that in both low and high cases we have
$\sum_{k=0}^M|F^k(\hat{T})|^{1+\xi} < C|I_0|^{1+\xi}$. We may
apply the usual method to show that this means that
$\sum_{k=1}^{n_i -m'} |f^{k+m'}(\hat{T})|^{1+\xi} <
C\sigma_{i,m}\max_{m'< k\le n_i} |f^k(T)|^\xi$.  So there is some
$C_{casc}$ \st
$$
\sum_{k=1}^{n_i- n_{i+m+1}}|f^{ k+n_{i+m+1}}(T)|^{1+\xi} <
C_{casc}\sigma_{i, m}\max_{n_{i+m}< k\le n_i} |f^k(T)|^\xi$$ as
required.
\end{proof}

\section{Exceptional case}

\label{sec:exceptional}

In the last section we dealt completely with the saddle node
cascade.  It is easily shown, for example applying
Lemma~\ref{lem:others} below to all branches, that following a
saddle node cascade we have a well bounded case, and so the
conclusions of Proposition~\ref{prop:entiresum} hold.  An
Ulam--Neumann cascade, however, is not always followed by a well
bounded case. We estimate the sum for $F_i$ in this alternative
case here. Most of the sum is dealt with using the methods for the
well bounded case, but we need some new techniques to deal with
two of the branches of $F_i$.

\begin{figure}[htp]
\begin{center}
\includegraphics[width=0.4\textwidth]{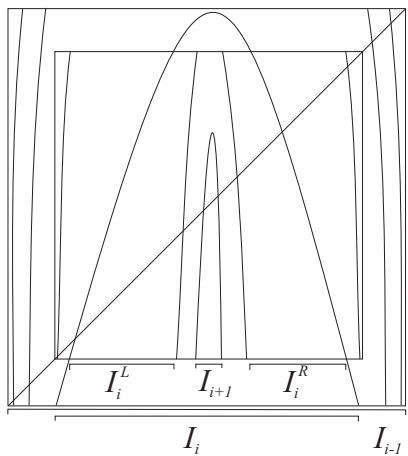}
\caption{The exceptional case.} \label{fig:exc}
\end{center}
\end{figure}

We consider the sum for $F_i$ where $F_{i-2}$ has a central return
and  $F_{i-1}$ has a high non--central return. The situation here
is only slightly different from the case considered in
Section~\ref{sec:bdd}, since we can prove that all domains of
$F_i$ are well inside $I_i$, except possibly two. Both of these
domains $I_i^j$ have $F_i|_{I_i^j} = F_{i-1}|_{I_i^j}$. We denote
the left--hand such interval by $I_i^L$ and the right--hand one by
$I_i^R$, see Figure~\ref{fig:exc}. These are the exceptional
domains.
If $I_{i-1}$ is a $\hat\chi$-scaled neighbourhood of $I_i$ then by
Theorem~\ref{thm:svarg} we know that $I_i$ is a
$\tilde{\hat\chi}$--scaled neighbourhood of both $I_i^L$ and
$I_i^R$, and we may proceed as in the well bounded case. But this
will not always be so if $I_{i-1}$ is at the end of a long
Ulam--Neumann cascade. So we will assume that $I_{i-1}$ is not a
$\tilde{\hat\chi}$-scaled neighbourhood of $I_{i}$. Without loss
of generality, we suppose that $F_{i-1}(c)$ is a maximum for
$F_{i-1}:I_i \to I_{i-1}$.

We are now ready to begin the proof of
Proposition~\ref{prop:entiresumex}.  The strategy for the proof is
as follows.

\begin{itemize}
\item Show there is some upper bound on $B(I_i, I_i^j)$ for $j
\neq L, R$.
\item State our main result in the proof:
Proposition~\ref{prop:boundI'}.  We suppose that we have some
interval $J \subset I_i^j$ for $j \neq L, R, 0$; $F_i(J), \ldots,
F_i^m(J) \subset I_i^L \cup I_i^R$; and $F_i^{m+1}(J) \subset
I_i^{j'}$ for $j' \neq L, R, 0$.  Then there exists some
$\lambda<1$ \st $B(I_i, J) < \lambda B(I_i, F_i^{m+1}(J))$.
Furthermore, $\sum_{k=1}^m|F_i^k(J)|< B(I_i, F_i^{m+1}(J))|I_i|$.
We are then able to prove Proposition~\ref{prop:entiresumex}. In
the rest of this section we prove Proposition~\ref{prop:boundI'}:
essentially we need an upper bound on $\sum_{k=1}^m|F_i^k(J)|$.

\item In Lemma~\ref{lem:dfbig} we show that
there exist an interval $V\subset I_i$ and $\gamma>1$ \st
$$|DF_i|_{(I_i^L \cup I_i^R) \setminus V} > \gamma.$$
This allows us to bound parts of the sum $\sum_{k=1}^m|F_i^k(J)|$
which lie in $(I_i^L \cup I_i^R) \sm V$.
\item We next focus on $V$.  We take first return maps to $V$ and
use decay of cross--ratios again to estimate sums of intervals in
$V$, see Lemma~\ref{lem:boundv}.  We can then complete the proof
of Proposition~\ref{prop:boundI'}
\end{itemize}

We first show in the following simple lemma that we have uniform
bounds on how deep the domains of $F_i$ are in $I_i$ for all
domains except $I_i^L, I_i^R$.

\begin{lem}
In the exceptional case outlined above, if $j\neq L, 0, R$ then
$I_i$ is a $\tilde{\hat{\chi}}$--scaled neighbourhood of $I_i^j$.
\label{lem:others}
\end{lem}

In fact, a similar result holds for the central domain too by
Theorem~\ref{thm:svarg}, but this is not important for us here.
This lemma proves that we can treat the case where $F_{i-2}$ is
central and $F_{i-1}$ is low and non--central as a well bounded
case.

As we shall see, the proof of this lemma is reminiscent of the
cascade case since we follow iterates of intervals along the
central branch of some $F_{i'}$.

\begin{proof}
There exists some maximal $i'<i$ \st $F_{i'-2}$ is non--central.
Then by Theorem~\ref{thm:svarg}, $I_{i'}$ is a $\hat\chi$--scaled
neighbourhood of $I_{i'+1}$.

For $j \neq L, R$ we will find $F_i|_{I_i^j}$ as a composition of
some branches of $F_{i'}$ in order to find some extensions.
$F_{i'}|_{I_{i'+1}}$ maps $I_i^j$ out of $I_i$ along the cascade,
through the sets $I_{i-1} \setminus I_i$, $I_{i-2} \setminus
I_{i-1}$ and so on, until it maps to some interval in $I_{i'+1}
\setminus I_{i'+2}$. Then this interval is mapped into some
$I_{i'}^{j'}$.  This then maps back into $I_{i'+1}$.  The process
may be repeated many times before $I_i^j$ is finally mapped back
to $I_i$.

So know that $F_i|_{I_i^j}$ is a composition of maps as follows.
Let $j_1 \neq 0$ be \st $(F_{i'}^{i-i'}|_{I_{i'+1}})(I_i^j)
\subset I_{i'}^{j_1}$.  Let $k_1= i-i'$.  If $F_i|_{I_i^j} =
(F_{i'}|_{I_{i'}^{j_1}}) (F_{i'}^{(i-i')}|_{I_{i'+1}} )|_{I_i^j}$
then we stop here; we say $r=1$.  Otherwise, let $k_2 \ge 0$ be
minimal \st $F_{i'}^{k_1+1+k_2}(I_i^j) \subset I_{i'} \setminus
I_{i'+1}$.  Let  $j_2 \neq 0$ be \st $F_{i'}^{ k_1+ 1+ k_2
}(I_i^j) \subset I_{i'}^{j_2}$.  If $F_i|_{I_i^j} = F_{i'}^{ k_1+1
+k_2+1}|_{I_i^j}$ then we stop here; we say $r=2$. Otherwise, we
continue this process until we finally return to $I_i$ and obtain
$k_r$.

Suppose that $r=1$.  That is,
$$F_i|_{I_i^j} = F_{i'}^{(i-i')+1}|_{I_i^j}.$$   Let $U$ denote
$F_{i'}^{(i-i')}(I_{i}^{j})$ and $U'$ denote $I_{i'}^{j_1}$.  Then
$F_{i'}(U) = I_i$ and $F_{i'}(U') = I_{i'}$. We know that $I_{i'}$
is a $\hat\chi$--scaled neighbourhood of $I_{i}$. So if we can
show that, taking the appropriate branch,
$(F_{i'}^{-(i-i')}|_{I_{i'+1}})(U') \subset I_i$, we know by
Theorem~\ref{thm:disjkobcr}(b) that $I_i$ is a
$\tilde{\hat{\chi}}$--scaled neighbourhood of $I_i^j$ (since all
the intervals we are concerned with are disjoint).  It is easy to
see that for this branch, $(F_{i'}^{-(i-i')}|_{I_{i'+1}})(U')
\subset I_i$ by the structure of the saddle node cascade since we
have $(F_{i'}^{-1}|_{I_{i'+1}})(U') \subset I_{i'+1} \setminus
I_{i'+2}$, $(F_{i'}^{-2}|_{I_{i'+1}})(U') \subset
I_{i'+2}\setminus I_{i'+3}$ and so on.  So the lemma is proved
when $r=1$.

In the more general case, where $r>1$ and $$F_i|_{I_i^j} =
F_{i'}^{\sum_{l=1}^r (k_l+1)}|_{I_i^j}$$ we may apply the same
idea, again using the disjointness of the domains of the first
return map, to prove that $I_i$ is a $\tilde{\hat{\chi}}$--scaled
neighbourhood of $I_i^j$.
\end{proof}
If necessary we adjust $\lambda$ so that
$\lambda(\tilde{\hat{\chi}})\le \lambda<1$.

By the above, if $I_i$ is a $\tilde{\hat\chi}$--scaled
neighbourhood of $I_i^L$ and $I_i^R$ then we can proceed with the
method in the well bounded case to prove
Proposition~\ref{prop:entiresumex}. But this is not generally the
case.  So for our work here, we may assume that $I_i$ is not a
$\tilde{\hat\chi}$--scaled neighbourhood of $I_i^L$ or $I_i^R$,
and that some iterate of $J$ enters $I_i^L \cup I_i^R$.

\begin{rmk}  In the previous sections we had uniform upper bounds
on the cross--ratio $B(I_i, I_i^j)$ for all $j$ and so we obtained
estimates on the decay of cross--ratios directly.  This was used
to estimate the sums of intervals.  The problem we often encounter
in this section is that sometimes we only get good estimates on
how cross--ratios decay and sometimes we only get good estimates
for the decay of the sizes of intervals.  But these estimates are
difficult to marry together directly, so we will have to split up
such cases.  The process is first described in the proof of
Proposition~\ref{prop:entiresumex} and again in the proof of
Lemma~\ref{lem:boundv}.  (As we will see later, this splitting
scheme deals with the cases where we enter $I_i^L \cup I_i^R$ from
$I_i$; $V$ from $I_i^L \cup I_i^R$; and $\Lambda$ from $V$.)
\end{rmk}
The principal result in this section is the following proposition.

\begin{figure}[htp]
\begin{center}
\includegraphics[width=0.8\textwidth]{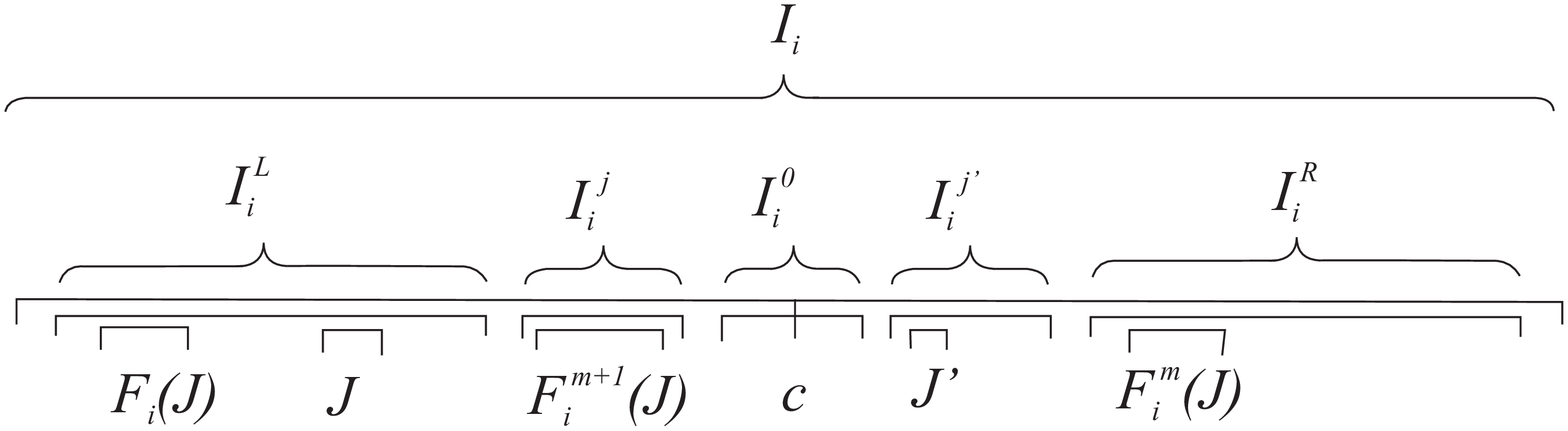}
\caption{An illustration of Proposition~\ref{prop:boundI'}.}
\label{fig:boundI'}
\end{center}
\end{figure}

\begin{prop}
If $J, F_i(J), \ldots, F_i^m(J) \subset I_i^L \cup I_i^R$ then
\begin{enumerate}
\item there exists some $0 \le \hat{m} < m$ \st
$\sum_{k=0}^m|F_i^k(J)| < C (|F_i^{m}(J)| + |F_i^{\hat{m}}(J)|)$;
\item  for some $\lambda<1$ independent of
$i$, if $F_i^{m+1}(J) \subset I_i^j$, $j \neq L,0, R$ then
\begin{enumerate}
\item $\sum_{k=0}^m|F_i^k(J)| < CB(I_i, F_i^{m+1}(J))|I_i|$;
\item letting $J'$ be the element of $F_i^{-1}(J)$ inside some
interval $I^{j'}$ for  $j' \neq L,0, R$ then we have $B(I_i, J') <
\lambda B(I_i, F_i^{m+2}(J'))$.
\end{enumerate}
\end{enumerate}
\label{prop:boundI'}
\end{prop}
See Figure~\ref{fig:boundI'} for a schematic representation of the
situation of this proposition. If necessary we will adjust the
$\lambda<1$ we use throughout this paper so that we may assume
that the proposition above holds for that $\lambda$.

\begin{proof}[Proof of Proposition~\ref{prop:entiresumex} assuming
Proposition~\ref{prop:boundI'}] As in the proof in the well
bounded case, we first show that we are principally concerned with
the intervals inside $I_i$. Again, the proof of this fact is a
slight modified version of the proof in the well bounded case.

Let $n_{i+1} < m_{1} <  \cdots < m_{j_i} = n_{i}$ be all the
integers between $n_{i+1}$ and $n_{i}$ \st $f^{m_{j}}(T) \subset
I_i \setminus I_{i+1}$ for $j=1, \ldots, j_i-1$ and let $m_{0} =
n_{i+1}$. Let $F_i:\bigcup_i U_i^j \to I_i$ be the first entry map
to $I_i$. As before, we will decompose the sum
$\sum_{i=n_{i+1}+1}^{n_i} |f^i(T)|$ as $\sum_{j=0}^{j_i-1}
\sum_{k=1}^{m_{j+1} - m_{j}} |f^{m_{j}+ k}(T)|$.

Suppose that $f^{m_{j}+1}(T) \subset U_i^j$ for some $U_i^j$.
Suppose further, that $F_i|_{U_i^j}=f^{i_j}$. Then there exists an
extension to $V_i^j \supset U_i^j$ so that $f^{i_j}: V_i^j \to
I_{i'-1}$ is a diffeomorphism, where $i'$ is defined in the proof
of Lemma~\ref{lem:others}. Then we have distortion bounds as
usual: $\frac{|f^k(f^{m_{j}+1}(T))|}{|f^k(U_i^j)|} \le C(\chi)
\frac{|f^{m_{j}+1}(T)|}{|I_i|}$.  Thus,
$$\sum_{k=1}^{m_{j+1} - m_{j}} |f^{m_{j}+ k}(T)|
< C(\chi) \sigma_{i} \frac{|f^{m_{j+1}}(T)|}{|I_{i}|}.$$

Therefore, $\sum_{j = n_{i+1}+1}^{n_i} |f^i(T)| < C(\chi)
\frac{\sigma_{i}}{|I_{i}|} \sum_{j=1}^{j_i}|f^{m_{i}}(T)|$. I.e.
we are principally interested in the sum
$\sum_{j=1}^{j_i}|f^{m_{i}}(T)|$, that is
$\sum_{k=0}^{j_i-1}|F_i^k(\hat{T})|$ where $\hat{T} =
f^{m_{1}}(T)$. In fact, we focus on bounding
$\sum_{k=0}^{j_i-2}|F_i^k(\hat{T})|$.

We split $\hat T, F_i(\hat T), \ldots, F_i^{j_i-2}(\hat{T})$ into
two groups: one for those intervals outside $I_i^L \cup I_i^R$ and
one for those inside $I_i^L \cup I_i^R$. Suppose that $J$ is an
interval \st for some $k \ge 0$, we have $F_i^{k}(J) \subset
I_i^{j}$ for some $j \neq L, 0, R$; then $F_i^{k+1}(J),
F_i^{k+2}(J), \ldots, F_i^{k'}(J) \subset I_i^L \cup I_i^R$ for
some $k'>k$; and finally $F_i^{k'+1}(J) \subset I_i^{j'}$ for some
$j' \neq L,0, R$. From the last part of
Proposition~\ref{prop:boundI'} we have
$$B(I_i, F_i^{k}(J))< \lambda B(I_i, F_i^{k'+1}(J)).$$ Therefore,
we can bound the sums of intervals which lie in the intervals
$I_i^j$ for all $j \neq L, R$ in a similar manner to that for the
well bounded case, independently of those intervals inside $I_i^L
\cup I_i^R$, as follows.

Given $k \ge 0$ \st $F_i^k(\hat{T}) \subset I_i^j$ for some $j
\neq L, 0, R$ we wish to estimate $|F_i^k(\hat{T})|$.  Let $0 \le
\hat{k} \le j_{i}-2$ be maximal \st $F_i^{\hat{k}}(\hat{T})
\subset I_i^{j'}$ for some $j' \neq L, R$.  Then we apply
Proposition~\ref{prop:boundI'} repeatedly to obtain $B(I_i,
F_i^{k}(\hat{T}))  < \lambda^l B(I_i, F_i^{\hat{k}}(\hat{T}))$ for
some $l \ge 0$.  The $l$ counts the number of times that
$F_i^{k+r}(\hat{T})$ lies outside $I_i^L \cup I_i^R$ for $0< r \le
\hat{k}$. Then
$$|F_i^k(\hat{T})| < \frac{|I_i|}{1+ \frac{2}{\lambda^l B(I_i,
F_i^{\hat{k}}(\hat{T}))}}.$$

We have two cases.  In the first case we have $\hat{k} = j_{i}-2$.
Then
\begin{eqnarray*}
B(I_i, F_i^{j_i-2}(\hat{T})) & < & B(I_i, I_i^{j'})
\frac{|F_i^{j_i-2}(\hat{T})|}{|I_i^{j'}|} <
\Delta(\tilde{\hat{\chi}})
\frac{|F_i^{j_i-2}(\hat{T})|}{|I_i^{j'}|} \\
& < & \Delta(\tilde{\hat{\chi}})C(\tilde{\hat{\chi}})
|F_i^{j_i-1}(\hat{T})|.
\end{eqnarray*}
Therefore, $|F_i^k(\hat{T})| <C\lambda^l |F_i^{j_i-1}(\hat{T})|$.
This suffices to prove an upper bound of the form
$C|F_i^{j_i-1}(\hat{T})|$ for the $F_i$--iterates of $\hat T$
outside $I_i^L \cup I_i^R$ in this case.

In the second case $\hat{k} < j_i-2$.  We have $$B(I_i,
F_i^{\hat{k}}(\hat{T})) < B(I_i,I_i^{j'})\frac{|
F_i^{\hat{k}}(\hat{T})|}{|I_i^{j'}|} < \frac{\Delta |
F_i^{\hat{k}}(\hat{T})|}{|I_i^{j'}|}.$$ Since
$|F_i^{\hat{k}}(\hat{T})| < C(\chi)| F_i^{\hat{k}+1}(\hat{T})|
\frac{|I_i^{j'}|}{|I_i|}$.  Therefore, in this case we have a
bound of the form $C|F_i^{\hat{k}+1}(\hat{T})|$ for the iterates
of $T$ outside $I_i^L \cup I_i^R$.

Finally we use the above information about sizes of intervals
outside $I_i^L \cup I_i^R$ to bound the sums of intervals inside
$I_i^L \cup I_i^R$ too.  In the first case above, we have a bound
of the form $C|F_i^{j_i-1}(\hat{T})|$ for the iterates of $T$ in
$I_i^L \cup I_i^R$.  In the second case above, we have a bound of
the form $C(|F_i^{\hat k}(\hat{T})| + |F_i^{\hat m}(\hat{T})|+
|F_i^{j_i-1}(\hat{T})|)$ for the iterates of $T$ in $I_i^L \cup
I_i^R$.

So in the worst case we have the bound $$C_{ex} \sigma_{i}
\left(\frac{|f^{n_i}(T)|}{|I_{i}|} + \frac{|f^{n_i, 2
}(T)|}{|I_{i}|}+\frac{|f^{n_i, 3}(T)|}{|I_{i}|} \right)$$ for the
sum $\sum_{k=n_{i+1}+1}^{n_i} |f^k(T)|$, as required.
\end{proof}

\subsection{Proof of Proposition~\ref{prop:boundI'}}

Denote the smallest interval containing both $I_i^L$ and $I_i^R$
by $I'_i$. Recall that we are assuming that the critical point is
a maximum for $F_{i-1}|_{I'_i}$.  (Recall that $F_i|_{I_i^L\cup
I_i^R} =F_{i-1}|_{I_i^L\cup I_i^R}$.)  This means that there is
some fixed point $p$ of $F_i$ in $I_i^R$. Clearly, there also
exists a point $p' \in I_i^L$ \st $F_i(p') =p$. Let $V:=(p', p)$.

We outline the proof of Proposition~\ref{prop:boundI'} as follows.
We suppose that some iterate of $J$ enters $V$.  Let $0 \le s_1
\le s_2 \le s_3$ be defined as follows. $F_i^{k}(J) \subset I'_i
\setminus V$ for $1 \le k \le s_1$; $F_i^{s_1+1}(J) \subset V \cap
(I_i^L \cup I_i^R)$; and $F_i^{s_2+k}(J) \subset  I'_i \setminus
V$ for $1 \le k \le s_3 -s_2$.  Any sum of the form $\sum_{k=0}^m
|F_i^k(J)|$ can be broken up into blocks consisting of such sums.

The scheme for proving Proposition~\ref{prop:boundI'} is to
firstly to show that $|DF_i|_{I'_i \setminus V}$ is uniformly
large. This is proved in Lemma~\ref{lem:dfbig} and helps to deal
with the sums $\sum_{k=0}^{s_1}|F_i^{k}(J)|$ and $\sum_{k=1}^{s_3-
s_2}|F_i^{s_2+ k}(J)|$. Then we have to prove that we have bounds
on the sums of intervals which return to $V$. This, proved in
Lemma~\ref{lem:boundv}, helps to deal with
$\sum_{k=1}^{s_2-s_1}|F_i^{s_1+ k}(J)|$.

Note that the proof of Proposition~\ref{prop:boundI'} is the only
time in this paper that we use the symmetry of the map (and it is
only a simplifying assumption).

\begin{lem}
There exists some $\gamma>1 $ independent of $i$ \st
$$|DF_i|_{I'_i \setminus V}> \gamma.$$
\label{lem:dfbig}
\end{lem}

\begin{proof}
We start by observing as in the last section that $|DF_i(p)|>
\rho_f$. By symmetry, $|DF_i(p')|>\rho_f$ too. Observe that
$I_i^L$ also contains a fixed point $q$ of $F_i$. We have
$|DF_i(q)|> \rho_f$ too. Furthermore, there exists a point $q' \in
I_i^R$ \st $F_i(q') = q$. From symmetry, $|DF_i(q')|> \rho_f$.

We can estimate $|DF_i|_{(p, q')}$ using the Minimum Principle as
follows. We use our $\mu$ given in \eqref{eq:mu} in place of
$\mu_g$. Then $|DF_i|_{(p, q')}> \mu^3\rho_f$. When $I_0$ is small
enough, $\mu$ is close to 1. Thus we may ensure that our intervals
are so small that $|DF_i|_{(p, q')}> \rho$ for some $\rho>1$. (To
fix precisely how small our intervals must be, we can, for
example, choose $\rho= \sqrt{\rho_f}$.) By symmetry, $|DF_i|_{(q,
p')}> \rho$.

We deal with the remaining part of the proof of the lemma by
showing that $F_i$ has large derivative when $x \in I'_i$ and
either $x<q$ or $x>q'$.  We use the following consequence of
Theorem~\ref{thm:svarg} and the Minimum Principle.

\begin{claim*}
There exists some $\gamma'=\gamma'(\chi)>1$ such that, denoting
$I_i^L= (l^-, l^+)$ and $I_i^R= (r^-, r^+)$, if $I_0$ is
sufficiently small and $B(I_i, I_i^L), B(I_i, I_i^R) $ are
sufficiently large then
$$|DF_i|_{(l^-,q)}, |DF_i|_{(q', r^+)}>\gamma'.$$
\end{claim*}

\begin{proof}
Let $\theta:=\frac
12\left(\frac{|I_{i'}|}{|I_{i'+1}|}-1\right)>\hat\chi$ where $i'$
is defined in the proof of Lemma~\ref{lem:others}.   We suppose
that $|DF_{i'}|_{I_{i'+1}\sm I_i} \le 1+2\theta$.  Then we prove
by induction that $\frac{|I_{i'+k}|}{|I_{i'+k+1}|} \ge 1+2\theta$
for $0 \le k <i-i'$.  By construction it is true for $k=0$.  We
assume that it is true for some $0\le k<i-i'-1$.  Then
\begin{align*}
\frac{|I_{i'+k+1}|}{|I_{i'+k+2}|} & \ge \frac{|I_{i'+k+2}|+
(\sup_{I_{i'+1}\sm I_{i}}|DF_{i'}|)^{-1} |I_{i'+k}\sm I_{i'+k+1}|
}{|I_{i'+k+2}|} \\
&\ge 1+\left(\frac{2\theta}{1+2\theta}\right)
\frac{|I_{i'+k+1}|}{|I_{i'+k+2}|}
\end{align*}
Then it is easy to see that $\frac{|I_{i'+k+1}|}{|I_{i'+k+2}|}\ge
1+2\theta$ as required.

In particular, we have proved that $|DF_{i'}|_{I_{i'+1}\sm I_i}
\le 1+2\theta$ implies that $I_i$ is a $\tilde\theta$--scaled
neighbourhood of both $I_i^L$ and $I_i^R$: a contradiction (since
$\tilde\theta>\tilde{\hat\chi}$). So there must exist some $x \in
I_{i'+1}\sm I_i$ \st $|DF_i(x)|\ge 1+2\theta > 1+2\hat\chi$.
Therefore, by Theorem~\ref{thm:min} and \eqref{eqn:sigma'sum} we
have
$$|DF_{i'}|_{(x_0, p)}
> \mu^3\min(1+2\hat\chi, \rho_f).$$
Choosing $|I_0|$ small we have some $\gamma'>1$ \st $|DF_i|_{(x_0,
q)} > \gamma'$.  In particular $|DF_i|_{(l^-, q)} > \gamma'$.
Similarly we can show $|DF_i|_{(q', r^+)} > \gamma'$.
\end{proof}
Letting $\gamma:= \min(\rho, \gamma')$, the lemma is proved.
\end{proof}

By the above, we will be able to estimate the sizes of iterates of
$T$ inside $(I_i^L \cup I_i^R) \setminus V$ as a  geometric sum.

We will need some real bounds for $V$. The following lemma, which
contrasts with Lemma~\ref{lem:dfbig}, will later be used to obtain
these bounds.

\begin{lem}
There exists some $\hat{C} = \hat{C}(\chi, |I'_i|)>0$, where
$\hat{C}(\chi, |I'_i|)$ tends to some constant $\hat{C}(\chi)$ as
$|I'_i| \to 0$, \st
$$|DF_i|_{I_i^L\cup I_i^R}< \hat{C}.$$
\label{lem:dfbound}
\end{lem}

\begin{proof} We work with $F_{i'}: I_{i'+1} \to I_{i'}$ where $i'$ is
defined in the proof of Lemma~\ref{lem:others}. There exists some
$m \ge 1$ \st $F_{i'}|_{I_{i'+1}} = f^m|_{I_{i'+1}}$. We can
decompose this map into two maps so that $F_{i'}= L\circ g$ where
$g = f|_{U_\phi}$, i.e $g(x) = f(c) - |\phi(x)|^\alpha$, and $L=
f^{m-1}: f(I_{i'+1}) \to I_{i'}$.

By Theorems~\ref{thm:disjkobcr}(a) and \ref{thm:svarg}(a) we have
$\frac{DL(x)}{DL(y)} < C(\chi)$ for $x, y \in f(I_{i+1})$. So
\begin{eqnarray*}
|DL(x)| & \le & C(\chi) \frac{|I_{i}|}{|f(I_{i+1})|} = C(\chi)
\frac{|I_{i}|}{\left|\phi\left(\frac{|I_{i+1}|}{2} \right)^\alpha
\right|}
\end{eqnarray*}
 for $x \in f(I_{i+1})$.  Also $$|Dg(x)| = \alpha |D\phi(x)|
|\phi(x)^{\alpha-1}| < \alpha \sup_{x \in
I_{i'+1}}|D\phi(x)|\left|\phi\left(\frac{|I_{i+1}|}{2}\right)\right|^{\alpha
-1}.$$ For $\hat{U} \subset U_\phi$ a small neighbourhood of $c$,
let ${\rm Dist}(\phi, \hat{U}):=\sup_{x, y \in
\hat{U}}\frac{|D\phi(x)|}{|D\phi(y)|}$. Observe that as $I'_i$
becomes smaller, ${\rm Dist}(\phi, I'_i)$ tends to 1. For $x \in
I_i^L\cup I_i^R$,
$$
|DF_i(x)|  <  \alpha C(\chi) \frac{\sup_{x \in I_{i+1}}|D\phi(x)|
|I_{i}|}{ \left|\phi\left(\frac{|I_{i+1}|}{2}\right)\right|} <
2\alpha C(\chi) {\rm Dist}(\phi, I'_i) \frac{ |I_{i}|}{
|I_{i+1}|}.$$ Since we have assumed that $\frac{ |I_{i}|}{
|I_{i+1}|}$ is bounded below, there is some constant $C>0$ \st for
all $x \in I_i'$,
$$|DF_i(x)| <  C C(\chi){\rm Dist}(\phi, I'_i).$$ Letting $\hat{C}(\chi,
|I'_i|) := CC(\chi) {\rm Dist}(\phi, I'_i)$ we have proved the
lemma.
\end{proof}

We denote the first return map to $V$ by $\hat{F_i}: \bigcup_j V^j
\to V$.  We first wish to find some control on the sizes of the
domains of $\hat{F_i}$.  Let $m_{V, j}$ be \st $\hat{F_i}|_{V^j} =
F_i^{m_{V, j}}|_{V^j}$.  The following lemma is key to proving
Proposition~\ref{prop:boundI'}.

\begin{lem}
If $F_i^{l_1}(J), \dots, F_i^{l_m}(J) \subset V\cap (I_i^L \cup
I_i^R)$ are all the iterates of $J$ up to $l_m$ which lie in
$V\cap (I_i^L \cup I_i^R)$, and all intermediate iterates
$F_i^k(J)$ for $k=0, 1, \ldots, l_m$ lie in $I_i^L \cup I_i^R$
then
$$\sum_{k=0}^{l_m}|F_i^k(J)| < C |F_i^{l_m}(J)|.$$
Furthermore, there exists $\lambda_V<1$ \st $|J| <C\lambda_V^{l_m
- m} |F_i^{l_m}(J)|$. \label{lem:boundv}
\end{lem}

\begin{proof} We split the sum as follows
$$\sum_{k=0}^{l_m}|F_i^k(J)| = \sum_{j=0}^{m-1} \sum_{k= 1}^{l_{j+1} - l_{j}}
|F_i^{l_{j} + k}(J)|$$ where we let $l_0=-1$. We know from
Lemma~\ref{lem:dfbig} that $|DF_i|_{I'_i \setminus V}> \gamma$ so
$$\sum_{k=1}^{l_{j+1} - l_{j}} |F_i^{l_{j}+ k}(J)| < |F_i^{l_{j+1}}(J)|
\sum_{k=0}^{l_{j+1} - l_{j}-1 } \gamma^{-k} <
\frac{|F_i^{l_{j+1}}(J)|}{1-\gamma^{-1}}.$$ Whence,
$$\sum_{k=0}^{l_m}|F_i^k(J)| <
\frac{1}{1-\gamma^{-1} } \sum_{j=0}^m |F_i^{l_j}(J)|.$$ So we only
need bound the sum of returns to $V$.

Denote the rightmost element of $\bigcup^jV^j$ by $V^1$ and the
leftmost element by $V^2$ (observe that $\hat{F_i}|_{V^1} =
F_i^2|_{V^1}$ and $\hat{F_i}|_{V^2} = F_i^2|_{V^2}$).  We get an
estimate on how deep each $V^j$ is inside $V$ for $j> 2$ because
$V^1$ and $V^2$ have some definite size compared to $|V|$; since
by Lemma~\ref{lem:dfbound} we know that $|V^1|, |V^2|
> \frac{|V|}{\hat{C}^2}$.  Therefore, there exists some
$\delta_0'$ depending only on $f$ \st $V$ is a $\delta_0'$--scaled
neighbourhood of $V^j$ for all $j>2$. So by
Lemma~\ref{lem:lambda'}, there exists some $\lambda_V'<1$
depending on $\delta_0'$ \st for any interval $J' \subset V^j$,
$B(V, J') < \lambda_V' B(V^j, J')$ for $j>2$ (in fact this is also
shown in the claim below). As usual we can use
Lemma~\ref{lem:lambda} to conclude that there exists some
$\lambda_V<1$ \st $B(V, J') < \lambda_V B(V, \hat{F_i}(J'))$. If
we remain away from $V^1$ and $V^2$, this fact and the usual
argument would be sufficient to obtain the required bound on sums.

We must deal with the case where iterates enter $V^1, \ V^2$.  The
idea is to split the situation into the case where intervals land
in a region where $|D\hat{F_i}|$ is large and the case when the
intervals land in a region where we don't have good estimates on
$|D\hat{F_i}|$.

We first focus on $V^2$.  We know from Theorem~\ref{thm:B} that
$|DF_i(p')|> \rho_f$ and so $|D\hat{F_i}(p')|> \rho_f^2$.  There
must also exist some fixed point $r$ of $\hat{F_i}$ in $V^2$ with
$|D\hat{F_i}(r)|> \rho_f$. Letting $\Lambda_2 := (p', r)$ and
applying the Minimum Principle as before, we obtain
$|D\hat{F_i}|_{\Lambda_2}> \rho$ for some $\rho>1$. Let $r'$ be
the point in $V^1$ \st $\hat{F_i}(r') = r$.  Then adjusting
$\rho>1$ if necessary, $|D\hat{F_i}|_{(r', p)}> \rho$. We define
$\Lambda_1$ to be the interval in $V^1$ which has
$\hat{F_i}(\Lambda_1) = V \setminus V^2$. Clearly $\Lambda_1
\subset (r', p)$, so $|D\hat{F_i}|_{\Lambda_1}> \rho$. For
convenience later, we let $\Lambda := \Lambda_1 \cup \Lambda_2$.

We are now ready to deal with bounding
$\sum_{k=0}^{m-1}|\hat{F_i}^k(J)|$.  Observe that
$\hat{F_i}^{m-1}(J)$ must be contained in some $V^j$.  Suppose
first that $j>2$; we deal with the case where $j =1$ or 2 later.
Suppose further that $J \subset V^{j'}$ and $j'>2$; here the other
case is similar.  We will again split up the sum. Let $N_0'=0$.
Let $N_1$ be minimal \st $\hat{F_i}^{N_1}(J) \cap \Lambda = \OOO$
and $\hat{F_i}^{N_1+1}(J) \subset \Lambda$.  Let $N_1'>N_1$ be
minimal \st $\hat{F_i}^{N_1'}(J) \subset \Lambda$ and
$\hat{F_i}^{N_1'+1}(J) \cap \Lambda = \OOO$.  In this way we
obtain $N_0'< N_1 < N_1' < \cdots < N_{M-1}< N_{M-1}'$ so that
\begin{eqnarray*}
\sum_{k=0}^{m-1}|\hat{F_i}^k(J)| & = & \sum_{j=0}^{M-1} \left(
\sum_{k=1}^{N_{j+1}- N_j'} |\hat{F_i}^{N_{j}'+k}(J)| +
\sum_{k=1}^{N_{j+1}'- N_{j+1}} |\hat{F_i}^{N_{j+1}+k}(J)| \right)
\\
& & + \sum_{k=1}^{N_{M}- N_{M-1}'} |\hat{F_i}^{N_{M-1}'+k}(J)|
\end{eqnarray*}
where $N_M =m-1$.  Observe that the first sum in the brackets
concerns intervals which land inside $\Lambda$ and the second sum
in the brackets concerns intervals in $V \setminus \Lambda$.  Then
$$\sum_{k=1}^{N_{j+1}'- N_{j+1}} |\hat{F_i}^{N_{j+1}+k}(J)| <
|\hat{F_i}^{N_{j+1}'}(J)|\sum_{k=0}^{N_{j+1} -N_{j+1}'-1 }
\rho^{-k} < \frac{C}{1-\rho^{-1}} |\hat{F_i}^{N_{j+1}'}(J)|
$$
for some $C$.

Now we consider $\sum_{k=1}^{N_{j+1}- N_j'}
|\hat{F_i}^{N_{j}'+k}(J)|$.  In fact we learn most from estimating
the sum $\sum_{k=1}^{N_{M}- N_{M-1}'}
|\hat{F_i}^{N_{M-1}'+k}(J)|$.  If necessary we make $\lambda_V<1$
smaller so that for $J \subset V^j \setminus \Lambda_j$ for
$j=1,2$ we have $B(V, J) < \lambda_V B(V, F_i(J))$. Then for $1
\le k< N_m - N_{M-1}'$,
$$B(V, \hat{F_i}^{N_{M-1}'+k}(J)) < \lambda_V^{N_{M}-N_{M-1}'-k} B(V,
\hat{F_i}^{N_M}(J)).$$  Recalling that $M = m-1$ we calculate
$B(V, \hat{F_i}^{m-1}(J)) < B(V, V^j)
\frac{|\hat{F_i}^{m-1}(J)|}{|V^j|}$. Letting $B_V := \max \left\{
\sup_{j>2}B(V, V^j), B(V, V^1 \setminus \Lambda_1), B(V, V^2
\setminus \Lambda_2)\right\} $, we obtain
$$|\hat{F_i}^{N_{M-1}'+k}(J)| < \frac{|V|}{1+ \frac{2 |V^j|}{
\lambda_V^{N_{M}-N_{M-1}'-k} B_V |\hat{F_i}^{m-1}(J)|}}.$$

Letting $\hat{B}_V := \frac{B_V}{B_V+2}$ we have
$$|\hat{F_i}^{N_{M-1}'+k}(J)| < \hat{B}_V \lambda_V^{N_{M}-N_{M-1}'+k}
\frac{|V|}{|V^j|} |\hat{F_i}^{m-1}(J)|.$$ Hence we have
$$\sum_{k=1}^{N_{j+1}- N_j'} |\hat{F_i}^{N_{j}'+k}(J)|< C
|\hat F_i^{m-1}(J)|.$$ We now estimate the other sums concerning
intervals outside $\Lambda$ as follows.  Let $\mu':=
\exp\left\{-\sigma'(I_0) \frac{|I_0|}{1-\rho^{-1}}\right\}$.
Suppose that $F_i^{N_{M-2}}(J) \subset V^j$.  Then taking the
appropriate branch, $\hat{F_i}^{N_{M-2} - N_{M-1}' -1}(V) \subset
V^j$ and
\begin{eqnarray*} B(V, \hat{F_i}^{N_{M-2}}(J)) & < & \lambda_V'
B(\hat{F_i}^{N_{M-2} - N_{M-1}' -1}(V), \hat{F_i}^{N_{M-2}}(J)) \\
& < & \frac{\lambda_V'}{\mu'} B(\hat{F_i}^{-1}(V), \hat{F_i}^{
N_{M-1}'}(J)) \\
& < & \frac{\lambda_V'}{\mu\mu'} B(V, \hat{F_i}^{ N_{M-1}'+1}(J))
\end{eqnarray*}
Shrinking $I_0$ if necessary, as usual, so that
$\frac{\lambda_V'}{\mu\mu'}=: \lambda_V<1$, we obtain
$$B(V, \hat{F_i}^{N_{M-2}}(J)) < \lambda_V B(V, \hat{F_i}^{
N_{M-1}'+1}(J)).$$ Clearly then we can proceed in bounding the sum
using the usual method of decaying cross--ratios.  So can bound
$\sum_{k=0}^{m-1}|\hat{F_i}^k(J)|$ above by $C|\hat F_i^{m-1}(J)|$
for this case.

To complete this case, we will bound $|\hat{F_i}^{m-1}(J)|$ in
terms of $|\hat{F_i}^{m}(J)|$.  We do this by constructing an
extension.  Let the left--hand and right--hand members of
$F_i^{-1}(p')$ be denoted by $b$ and $b'$ respectively. Denote
$(b, b')$ by $V'$.  By Lemma~\ref{lem:dfbound}, $V'$ is a
$\delta_{V'}$--scaled neighbourhood of $V$ where $\delta_{V'}$
depends only on $f$.

\begin{claim*}
For all domains $V^j$, $j>2$ there exists an extension to some
interval $U^j \supset V^j$ \st $U^j \subset V$ and $F_i^{m_{V,
j}}: U^j \to V'$ is a diffeomorphism. \label{claim:vext}
\end{claim*}

\begin{proof}
For $j>2$ the return maps are a composition of $F_i|_V$ followed
by $F_i|_{I_i^R}$ and then some number of iterates of
$F_i|_{I_i^L}$. So $\hat{F_i}^{-1}$ must pull $V'$ back into
$I_i^L$. Observe that this element of $F_i^{-1}(V')$ is below $p'$
(and clearly away from $F_i(c)$). Any further pullbacks in $I_i^L$
remain below $p'$ also. Therefore when some element $F_i^{-k}(V')$
is finally pulled back into $I_i^R$, it is mapped above $p$ and
remains away from $F_i(c)$. Therefore we have elements of
$F_i^{-k-2}(V')$ mapping inside $V$ which don't contain $c$.
\end{proof}

By the above claim and Theorem~\ref{thm:disjkobcr}  we have some
$C>0$ depending only on $f$ \st if $j>2$, $$\frac{1}{C}
\frac{|V|}{|V^j|} \le |D\hat{F_i}|_{V^j} \le C
\frac{|V|}{|V^j|}.$$ (Recall that we are assuming that
$F_i^{m-1}(V) \cap \Lambda = \OOO$.)

Therefore,
$$\sum_{k=1}^{N_{M}- N_{M-1}'} |\hat{F_i}^{N_{M-1}'+k}(J)| < C
|\hat{F_i}^{m}(J)|.$$

There remains a further case to consider. Above we assumed
$\hat{F_i}^{m-1}(J) \subset V^j$ where $j>2$.  But if $j \in \{1,
2\}$ we have two cases. We first note that if $ F_i^{l_m}(J) \cap
\{r, r'\} = \OOO$ then the intervals we are concerned with are
either completely inside $\Lambda_2, \Lambda_1$ or completely
inside $V \setminus (\Lambda_2 \cup \Lambda_1)$.  Then we may
proceed as above.  But if $F_i^k(J)$ contains $r$ or $r'$ then we
split $F_i^k(J)$ into two intervals, with this periodic point at
their intersection.  We may then apply the procedure above to
estimate the size of each interval.  We need only apply this
splitting argument once since if we intersect a periodic point of
$\hat{F_i}$ once, we must stay there for all time under iteration
by $\hat{F_i}$.  Thus we need only alter our constants by a factor
of 2 to deal with this case.  Note that we only have one sum where
this problem could occur:
$\sum_{k=1}^{N_M'-N_M}|\hat{F_i}^{N_M+k}(J)|$ where $N_M'=m$. This
is because $r$ is a fixed point for $\hat{F_i}$.

Clearly, we can use the cross--ratio argument as usual to obtain
the estimate $|F_i^{l_1}(J)|< \lambda_V^{m-1}C|F_i^{l_m}(J)|$, so
$|J|< \lambda_V^{m-1}C|F_i^{l_m}(J)|$. \end{proof}

We may adjust our usual $\lambda$ so that $\lambda_V \le \lambda
<1$.

\begin{proof}[Proof of Proposition~\ref{prop:boundI'}] Suppose first
that $F_i^{m+1}(J) \subset I_i^j$ for $j \neq L, R$.  Then, in
particular, we can be sure that $F_i^m(J)$ does not contain $p$ or
$p'$.  Then we also know that none of $F_i^k(J)$  contain $p$ or
$p'$ for $0 \le k \le m-1$. This means that we can be sure that
all the intervals we consider are either contained in $V$ or are
disjoint from $V$.

Recall that $0 \le s_1 < s_2 \le s_3 = m$ are defined as follows.
(We suppose that some iterate of $J$ enters $V$: otherwise the
proof is simpler.) $F_i^{k}(J) \subset I'_i \setminus V$ for $1
\le k \le s_1$; $F_i^{s_1+1}(J) \subset V \cap (I_i^L \cup
I_i^R)$; and $F_i^{s_2}(J) \subset V \cap (I_i^L \cup I_i^R)$,
$F_i^{s_2+k}(J) \subset I'_i \setminus V$ for $1 \le k \le s_3
-s_2$.

Then if $s_3 > s_2$,
$$\sum_{k=1}^{s_3-s_2}|F_i^{s_2+ k}(J)| < |F_i^{s_3}(J)|
\sum_{k=0}^{s_3-s_2-1} \gamma^{-k} < C|F_i^{s_3}(J)|,$$ by
Lemma~\ref{lem:dfbig}.

From Lemma~\ref{lem:boundv},
$$\sum_{k=1}^{s_2-s_1}|F_i^{s_1+ k}(J)|< C|F_i^{s_2}(J)|$$ and
$|F_i^{s_1+1}(J)|< C |F_i^{s_2}(J)|$.

Also \begin{eqnarray*} \sum_{k=0}^{s_1}|F_i^{k}(J)| & < &
\gamma^{-1}|F_i^{s_1+1}(J)| \sum_{k=0}^{s_1} \gamma^{-k} <
C|F_i^{s_2}(J)|.
\end{eqnarray*}
Therefore, $$\sum_{k=0}^{s_2}|F_i^{k}(J)| < C|F_i^{s_2}(J)|.$$ If
$s_3>s_2$ then
$$\sum_{k=0}^{s_3}|F_i^{ k}(J)| < C\left(|F_i^{s_3}(J)|+
|F_i^{s_2}(J)|\right).$$  Therefore, the first part of the
proposition is proved.

Now if $F_i^{m+1}(J) \subset I_i^j$ for $j \neq L, R, 0$ then
recalling that $s_3 = m$ we will obtain an estimate for
$|F_i^{s_2}(J)|$ in terms of $B(I_i, F_i^{m+1}(J))$.  $$B(I_i,
F_i^{s_2}(J)) < B(F_i^{-s_3+s_2}(I_i), F_i^{s_2}(J)) <
\frac{B(I_i, F_i^{m}(J))}{\mu} < \frac{B(I_i,
F_i^{m+1}(J))}{\mu^2}.$$ We are allowed to use $\mu$ here since
all intermediate intervals must be disjoint (otherwise we would
have to pass through $V$ again). Therefore $$|F_i^{s_2}(J)|<
\frac{|I_i|}{1+\frac{2\mu^2}{B(I_i, F_i^{m+1}(J))}} < C|I_i|
B(I_i, F_i^{m+1}(J)).$$  Similarly we can show that $|F_i^{m}(J)|<
C|I_i| B(I_i, F_i^{m+1}(J))$.  Therefore
$$\sum_{k=0}^{s_3}|F_i^{k}(J)| < C|I_i| B(I_i, F_i^{m+1}(J))<
C_1|I_i|$$ for some $C_1>0$.

We now prove the final part of the proposition.  Clearly for any
run of intervals $F_i(J), \ldots, F_i^k(J) \subset I_i^L \cup
I_i^R$, considering the branch of $F_i^{-k}$ which follows the
iterates of $J$, we have $B(F_i^k, F_i^{-k}(I_i), J)
> \mu''$ where $\mu'':=
\exp\{-C_1\sigma'(|I_0|)|I_0| \}$. We consider the branch of
$F_i^{-m-2}$ which follows the backward orbit of $F_i^{m+1}(J)$.
Clearly, $F_i^{-m-2}(I_i)$ is strictly inside $I_i^{j}$. Thus,
\begin{eqnarray*}
B(I_i, J') & < & \lambda' B(I_i^{j}, J') < \lambda'
B(F_i^{-m-2}(I_i), J') < \frac{\lambda'}{\mu''} B(F_i^{-1}(I_i),
F_i^{m+1}(J'))  \\ & < & \frac{\lambda'}{\mu'' \mu} B(I_i,
F_i^{m+2}(J')).
\end{eqnarray*}
For $|I_0|$ small enough, we can alter the usual $\lambda$
slightly so that $\frac{\lambda'}{\mu'' \mu} \le \lambda$ and
still ensure that $\lambda<1$.  Thus, $B(I_i, J')< \lambda B(I_i,
F_i^{m+2}(J'))$ as required.

When we do not escape $I_i^L \cup I_i^R$ then we may have some
intersection with $p$ or $p'$.  In this case, we split our
interval in two and estimate the size of each piece as above. We
need only apply this idea once, so we can change our constants to
cater for this case too.  In this case, part 2 of the proposition
doesn't occur.
\end{proof}

\section{Proof of the main theorem in the non--infinitely
renormalisable case} \label{sec:proof of cr non-inf}

We recall that $B(f^n, T, J)> \exp\{-C
\sum_{k=0}^{n-1}|f^k(T)|^{1+\eta}\}$ when $f \in C^{2+\eta}$. We
will find a bound on the sum $\sum_{k=0}^{n-1}|f^k(T)|^{1+\eta}$
by using the main propositions above and also finding some decay
property for the size of the domains of $F_i$ for some values of
$i$. We assume that $f^k(T)\cap\bd I_j\neq \OOO$ only within a
cascade case (i.e. when there exist $i,m$ \st $F_i$ is in a
cascade case and $f^k(T) \subset I_i\sm I_{i+m}$).  It is easy to
see how to extend the proof when this is not true.

Let $F_i : \bigcup_j U_i^j \to I_i$ be the first entry map to
$I_i$ (we include the branches of the first return map in this
case too). For $i<j$ and an interval $V$, we define $S(i, j, V)$
to be the maximum of $|f^{i+1}(V)|, |f^{i+2}(V)|, \ldots,
|f^{j}(V)|$. We will consider $S(n_{i+1}, n_i, T)$.  Let $n(i,j)$
be \st $F_i|_{U_i^j}= f^{n(i,j)}|_{U_i^j}$. Now let $U_{i}^{s(i)}$
be the interval for which $S\left(0, n(i,j), U_{i}^{j}\right)$ is
maximal. Let $\hat n(i) = n(i,s(i))$. Clearly,
$$S(n_{i+1}, n_i, T) \le S\left(0, \hat n(i), U_{i}^{s(i)}\right).$$
We would like to show that for certain $i$, this quantity decays
with $i$ in a controlled way.

We start by assuming that $F_{i-1}$ is in a well bounded case. We
have two cases. Firstly, suppose that $U_{i}^{s(i)} \subset I_i$.
Then since $F_{i-1}$ is in a well bounded case, we have
$|U_{i}^{s(i)}| < \frac{|I_{i-1}|}{1+2\chi}$.  Since $I_i$ is a
domain of the first return map to $I_{i-1}$ we have
$$|U_i^{s(i)}| < \frac{S\left(0, \hat n(i-1),
U_{i-1}^{s(i-1)}\right)}{1+2\chi}.$$

Now assume that $U_{i}^{s(i)} \cap I_i = \OOO$. Then there exists
some extension $V_i \supset U_{i}^{s(i)}$ \st $f^{n(s(i))}:V_i \to
I_{i-1}$ is a diffeomorphism.  We will show that $U_i^{s(i)}$ is
uniformly smaller than $V_i$.  By \eqref{eqn:sigma'sum} we know
that $B(V_i, U_i^{s(i)}) < \frac{B(I_{i-1}, I_i)}{\mu}$ for $\mu$
as in \eqref{eq:mu}. Thus, by Lemma~\ref{lem:Delta}, $|U_i^{s(i)}|
< \frac{|V_i|}{1+\frac{2\mu}{\Delta(\chi)}}$.  Since $V_i$ is a
first return domain to $I_{i-1}$ we have
$$|U_i^{s(i)}| < \frac{S\left(0, \hat n(i-1),
U_{i-1}^{s(i-1)}\right)}{1+\frac{2\mu}{\Delta(\chi)}}.$$

Let $\gamma := \max\left(\frac{1}{1+2\chi},
\frac{1}{1+\frac{2\mu}{\Delta(\chi)}} \right)$. Clearly
$\gamma<1$. So $$S\left(0, \hat n(i), U_i^{s(i)}\right) < \gamma
S\left(0, \hat n(i-1), U_{i-1}^{s(i-1)}\right).$$

We let $C_{all} = \max(C_{wb}, C_{casc}, 3C_{ex})$. Note that by
disjointness, all $\sigma_i, \sigma_{i,m}<1$.  If $f \in
NF^{2+\eta}$ and $F_{i-1}$ is well bounded, we have
\begin{align*} & B(f^{n_i-n_{i+1}}, f^{n_{i+1}+1}(T),
f^{n_{i+1}+1}(J))  \\
&  \hspace{4cm} \ge \exp\left\{-C\left(S( n_{i+1}, n_i,
T)\right)^\eta
\sum_{k=1}^{n_i-n_{i+1}} |f^{k+n_{i+1}}(T)| \right\} \\
& \hspace{4cm} > \exp\left\{-C\left(S\left(0, \hat n(i),
U_i^{s\left(i\right)}\right)\right)^\eta
C_{all} \right\} \\
& \hspace{4cm} > \exp\left\{-C\left(\gamma S\left(0, \hat n(i-1),
U_{i-1}^{s\left(i-1\right)}\right)\right)^\eta C_{all} \right\}.
\end{align*}

If we are not in the infinite cascade case then the sums for $F_i,
F_{i+1}, \ldots$ can be broken into blocks consisting of a
cascade; possibly followed by an exceptional case; followed by one
or more well bounded cases. So suppose that $F_i$ is well bounded,
$F_{i}, F_{i+1}, \ldots, F_{i+m-1}$ have central returns,
$F_{i+m}$ has a non--central return and $F_{i+m+1}$ is an
exceptional case. So note that, in particular, $F_{i+m+2}$ must be
well bounded.  Then,
\begin{eqnarray*} S\left(0, \hat n(i+m+3),
 U_{i+m+3}^{s(i+m+3)}\right) <
\gamma S\left(0, \hat n(i+m+2), U_{i+m+2}^{s(i+m+2)}\right),
\ldots \hspace{1cm} \\
\ldots, \gamma S\left(0, \hat n(i+1), U_{i+1}^{s(i+1)}\right) <
\gamma^2 S\left(0, \hat n(i), U_{i}^{s(i)}\right).
\end{eqnarray*}
Therefore, we have
\begin{eqnarray*} B(f^n, T, J) & > &
\exp\left\{-C\sum_{k=0}^{n-1}|f^k(T)|^{1+\eta} \right\} \\
&>& \exp\left\{-CC_{all}\left( S\left(0, \hat n(0),
U_{0}^{s(0)}\right) \right)^\eta
\sum_{k=0}^\infty \gamma^{k\eta} \right\} \\
&>& \exp\left\{-CC_{all}
\frac{\left(\sigma'(|I_0|)\right)^\eta}{1-\gamma^\eta}\right\}.
\end{eqnarray*}
Hence it is easy to see that for any $0<K<1$, if $I_0$ is the
central domain of a first return map to some $I_{-1}$, $I_0$ is
sufficiently small and $F_{-1}$ is non--central, then we may bound
$B(f^n, T, J)$ below by $K$.

Note that we can always start with a well--bounded case when we
don't have an infinite cascade.  We simply induce on a nice
interval finitely many times until we obtain a non--central return
and thus obtain a suitable $I_{-1}$. We consider the infinite
cascade case in the next section.

The second part of Theorem~\ref{thm:cr}, concerning $A(f^n, T,
J)$, is proved in the same way.

\section{Infinite cascade case}
\label{sec:infcasc}

\begin{figure}[htp]
\begin{center}
\includegraphics[height=0.35\textheight]{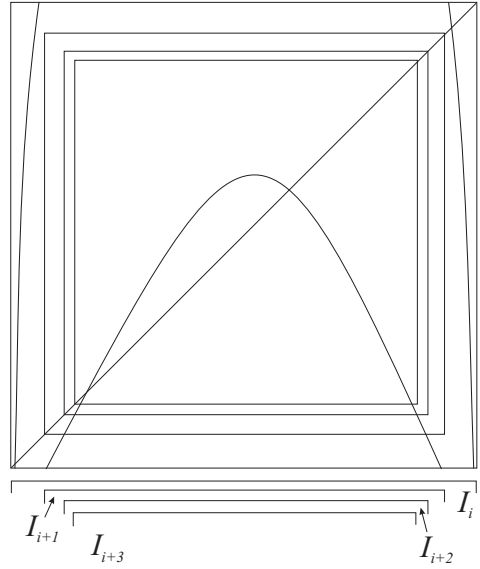}
\caption{An infinite cascade.} \label{fig:infcasc}
\end{center}
\end{figure}

Here we consider the case where we have some $I_0$ \st $F_i$ are
central for $i=0,1, \ldots$.  In this case we will find that
$\frac{|I_{i+1}|}{|I_i|}$ gets very close to 1. See
Figure~\ref{fig:infcasc} for an example of such a map. In
particular, $I_i$ will not shrink down to a point (the critical
point c) as $i$ increases so we can't use the method above to
bound sums of intervals which land very close to $c$. The
principal tool here is an extension given by a result of
\cite{kozsch}.  We will not supply all the details of our proof of
Theorem~\ref{thm:cr} in this case since the techniques are mostly
the same as applied in the previous sections.
We start by letting $I_0$ be any nice interval about $c$.  We
assume that we have some infinite cascade.  This means that for a
nice interval $I_0 \ni c$, $F_i$ is central (and high) for all
$i$, where $F_i$ is defined in the usual way. The main idea here
is that we can still find good bounds on some interval $I_{0, 0}$
and then apply the methods of Section~\ref{sec:cascade} to it.
Then we need to find another interval $I_{1, 0}$ around $c$ which
is smaller than all $I_{0, i}$, also has good bounds and is
uniformly smaller than $I_{0, 0}$. In such a way, we obtain a
sequence of intervals $I_{i, 0}$ which can each  be treated as in
the high cascade case above, and which shrink uniformly to the
critical point. Clearly $F_{i,j}$ will always be central and high
for all $i, j \ge 0$.

\begin{prop}
For $f \in NF^2$, and $\xi>0$ there exists some $C_{inf}>0$ \st
for any small $I_{0,0}$ defined as above, $T \subset I_{0,0}$
implies
$$\sum_{k=0}^{n-1}|f^k(T)|^{1+\xi}<C_{inf}.$$
\label{prop:infsum}
\end{prop}

Clearly this completes the proof of Theorem~\ref{thm:cr} in this
case.

\begin{proof}
We will prove this with a series of lemmas.

For all $i$ the central branch of $F_i$ has two fixed points,
$q_0$ and $p_0$ to the left and right of $c$ respectively (as
usual, we assume that $F_i(c)$ is a maximum for $F_i|_{I_{i+1}}$).
We let $q_0'$ be the point in $I_{i+1}$ not equal to $q_0$ which
maps by $F_i$ to $q_0$.  We define $p_0'$ similarly.  We define
$I_{0,0}$ to be $(p_0',p_0)$.  Let $F_{0, 0}:\bigcup_{j} I_{0,
0}^j \to I_{0, 0}$ be the first return map to $I_{0,0}$ (where
$I_{0, 0}^0$ is the central domain). We have the following lemma.

\begin{lem}
There exists some $\hat{\chi}>0$ depending only on $f$ \st
$I_{0,0}$ is a $\hat{\chi}$--scaled neighbourhood of every domain
$I_{0,0}^j$ which has $\bd I_{0,0}^j \cap \bd I_{0,0} = \OOO$.
\label{lem:hatchi}
\end{lem}

\begin{proof}
Clearly, $I_i$ tends to $(q_0, q_0')$.  So we denote $(q_0, q_0')$
by $I_{\infty}$. We will first show that $I_\infty$ is uniformly
larger than $I_{0,0}$ and then show that all except two
non--central domains of the first entry map to $I_{0,0}$ have an
extension to $I_\infty$ and show what this means for $I_{0, 0}^j$.
These two domains are the ones with either $p_0$ or $p_0'$ in the
closure.

In a similar manner to the exceptional case, we will find an upper
bound for $|DF_i|_{I_{i+1}}$. This will allow us to get good
bounds for the first return map to $I_{0,0}$

For large $i$, the ratio $I_i$ has $\frac{|I_{i+1}|}{|I_i|}$ close
to 1. The following lemma, an adaptation of Lemma 7.2 of
\cite{kozsch}, allows us to bound $|DF_i|_{I_{i+1}}$.

\begin{lem}
If $f \in NF^2$ then there exist constants $0< \tau_2<1$ and
$\tau_3>0$ with the following property.  If $T$ is any
sufficiently small nice interval around the critical point, $R_T$
is the first entry map to $T$ and its central domain $J$ is
sufficiently big, i.e. $\frac{|J|}{|T|}>\tau_2$, then there is an
interval $W$ which is a $\tau_3$--scaled neighbourhood of the
interval $T$ \st if $c \in R_T(J)$ then the range of any branch of
$R_T:V \to T$ can be extended to $W$ provided that $V$ is not $J$.
\label{lem:bigcentral}
\end{lem}

This lemma is only needed as a $C^3$ result in \cite{kozsch}, but
it easily extends to our $C^2$ case.

It is straightforward to see that the above lemma is sufficient to
prove a version of Lemma~\ref{lem:dfbound} in our case.  That is,
for large $i$, there exists some $\hat{C}'$ \st
$|DF_i|_{I_{i+1}}<\hat{C}'$.  This implies that there exists some
$0<\theta<1$ depending only on $f$ \st $|I_{0,0}|<
\theta|I_\infty|$ and, equivalently, some $\delta>0$ \st
$I_\infty$ is a $\delta$--scaled neighbourhood of $I_{0,0}$.

Now, for the moment we let $F_{0, 0}$ also denote the first entry
map and $\bigcup_j I_{0,0}^j$ also include the first entry
domains.  We will show that many of the branches have an extension
to a uniformly larger domain. Suppose that there exists a domain
$I_{0,0}^j$ with $\overline{I_{0,0}} \cap \overline{I_{0,0}^j} =
\OOO$ \st $F_{0,0}:I_{0,0}^j \to I_{0,0}$ does not have an
extension to $I_\infty$. That is, supposing $F_{0,0}|_{I_{0,0}^j}
= f^{n(j)}|_{I_{0,0}^j}$, there is no interval $V\supset
I_{0,0}^j$ \st $f^{n(j)}:V \to I_\infty$ is a diffeomorphism. Let
$0 \le k \le n(j)-1$ be maximal \st $f^{n(j)-k}: f^k(I_{0,0}^j)
\to I_{0,0}$ has no extension to $I_\infty$. Clearly if $I_{0, 0}$
is small $f: f^{n(j)-1}(I_{0,0}^j) \to I_{0,0}$ always has an
extension, so $k < n(j)-1$.  Then there exists some interval $W
\supset f^{k+1}(I_{0,0}^j)$ \st $f^{n(j)-k-1}:W \to I_\infty$ is a
diffeomorphism and the element $W'$ of $f^{-1}(V)$ containing
$f^{k}(I_{0,0}^j)$ contains $c$.

Since $I_\infty$ is a nice interval, $W'\subset I_\infty$.  We
also know that $f^{k}(I_{0,0}^j) \subset I_\infty \setminus
I_{0,0}$. Therefore $W'$ contains  either $p_0$ or $p_0'$.  But
then either $f^{n(j)-k-1}(p_0)$ or $f^{n(j)-k-1}(p_0')$ is
contained in $I_\infty \setminus I_{0,0}$ which is not possible.

Consider $I_{0,0}^j$ for some $j \neq 0$ where $I_{0,0}^j \subset
I_{0,0}$ is a domain of the first return map.  We will show that
this domain is uniformly deep inside $I_{0,0}$.  There exists some
$V \supset f(I_{0,0}^j)$, where $f^{n(j)}:V \to I_\infty$ is a
diffeomorphism and  $V$ is a $\tilde\delta$--scaled neighbourhood
of $f(I_{0,0}^j)$.  Let $V'$ be the maximal interval around
$I_{0,0}^j$ \st $f(V')=V$.  We show that $V' \subset I_{0,0}$. Let
$V(f(c))$ denote the maximal interval around $f(c)$ which pulls
back by $f^{-1}$ to $I_{0,0}$. If $V$ is not contained in
$V(f(c))$ then either $p_0$ or $p_0'$ is contained in $V'$. Thus,
$f^{n(j)}(p_0)$ or $f^{n(j)}(p_0')$ lies in $I_\infty \setminus
I_{0,0}$, a contradiction. So $V' \subset I_{0,0}$ and $I_{0,0}$
is a $\delta'$--scaled neighbourhood of $I_{0,0}^j$ where $\delta'
= \min\left(\tilde\delta, \frac{1}{2}\right)$. The case of the
central branch follows in the usual manner.
\end{proof}

So we are in a type of high cascade case for $F_{0,0}$.  Note that
the branches with $p_0$ or $p_0'$ in their closure can be dealt
with in the same way as the domains $V^1,V^2$ were dealt with in
the exceptional case.

We may assume that $F_{0,0}$ has an infinite cascade and is high
too. Let $F_{0,1}$ be the first return map to $I_{0,0}$ and so on,
so we obtain $I_{0,i}$. We sum for $F_{0,0}, F_{0,1}, \ldots$ as
in the high cascade case. We let $q_1, q_1', p_1, p_1'$ be defined
as above for the fixed points of $F_{0,0}|_{I_{0,1}}$. We let
$I_{0, \infty}$ denote $(q_1, q_1')$. We may apply the same ideas
as above to find some new interval $I_{1, 0}:= (p_1, p_1')$ which
has $|I_{1,0}|< \theta |I_{0, \infty}|$.  We may define $I_{i,j}$
for $i \ge 2$, and $0 \le j \le \infty$ in a similar way.

Let $f^{N_i}(T)$ be the last iterate of $T$ which lies inside
$I_{i,0}$.  Let $N_i'$ be the maximal integer $N_i \ge
N_i'>N_{i+1}$ \st $f^{N_i'}(T)$ is not in $I_{i,0} \sm I_{i,
\infty}$. Then these arguments prove the following lemma.

\begin{lem}
There exists some $C>0$ \st
$$\sum_{k=1}^{N_i-N_i'}|f^{k+N_i'}(T)| < C\hat{\sigma}_i$$
where $\hat{\sigma}_i$ is defined as follows.  Let $\sigma_i :=
\sup_{V \in {\rm dom} F_{i,0}} \sum_{j=1}^{n(V)} |f^j(V)|$ (and
$n(V)$ is defined as $k$ where $F_{i,0}|_V = f^k$).  Let $\hat{V}
\subset I_{i, 0} \setminus I_{i, 1}$  be an interval \st
$f^{\hat{n}}(\hat{V})$ is one of the connected components of
$I_{i, 0} \setminus I_{i, 1}$ and $f^j(\hat{V})$ is disjoint from
both $I_{i, 0} \setminus I_{i, 1}$  and $I_{i+1, 0}$ for $0< j<
\hat{n}(\hat{V})$. Then $\hat{\sigma}_{i}$ is the supremum of all
such sums $\sum_{j=1}^{\hat{n}(\hat{V})}|f^j(\hat{V})|$ and
$\sigma_i$. \label{lem:infsum}
\end{lem}

Now we consider $\sum_{k=1}^{N_i'-N_{i+1}}|f^{k+N_{i+1}}(T)|$.  If
none of these intervals contain $p_i,q_i$ then we are in
$I_{i,\infty} \sm I_{i+1,0}$.  By the Minimum principle,
$|DF_i|_{I_{i , \infty} \setminus I_{ i+1, 0}}$ is uniformly
greater than 1.  So we can easily bound our sum.  If none of our
intervals contains $p_i$, but some $f^{k+N_{i+1}}(T)$ contains
$q_0, q_0'$ we can split $f^{k+N_{i+1}}(T)$ at $q_0$ or $q_0'$
into two intervals. It is easy to see that there is some $C>0$ \st
$\sum_{k=1}^{N_i'-N_{i+1}}|f^{k+N_{i+1}}(T)|< C\hat\sigma_i$.  If
$p_0,p_0'$ is contained in some $|f^{k+N_{i+1}}(T)|$ then we must
split the interval at $p_0$ or $p_0'$.  Note that we may have to
split the interval $|f^{k+N_{i+1}}(T)|$ at arbitrarily many $p_i,
q_i'$ or $p_i', q_i$.  Therefore,
$\sum_{k=1}^{N_i'-N_{i+1}}|f^{k+N_{i+1}}(T)|^{1+\xi}<
C\sum_{k=i}^\infty (k-i)S(N_k, \hat n (k), T)^\xi\hat\sigma_k$
where $S$ and $\hat n$ are defined analogously to
Section~\ref{sec:proof of cr non-inf}.   As before there is some
constant $0<\theta'<1$, here depending on $\theta$ rather than
$\gamma$ such that $\theta'$ governs the decay of $S(N_i, \hat
n(i), T)$.  Hence, we can put this estimate together with
Lemma~\ref{lem:infsum} to get
$$\sum_{k=1}^{N_i-N_{i+1}}|f^{k+N_{i+1}}(T)|^{1+\xi}
< C S(N_i, \hat n (i), T)^\xi \sum_{k=0}^\infty k\theta'^{k\xi}.$$
Similarly to before, we can conclude that there exists some
$C_{inf}>0$ \st
$$\sum_{k=0}^{n} |f^{k}(T)|^{1+\xi} < C_{inf}.$$
\end{proof}

\appendix

\section{Proof of the Yoccoz Lemma}

\label{sec:yoc} We recall the lemma.

{\bf Lemma~\ref{lem:yoc}} {\it Suppose that $f \in NF^2$.  Then
for all $\delta, \delta'>0$ there exists $C>0$ \st if $I_0$ is a
nice interval \st \begin{enumerate}
\item $I_0$ is a $\delta$--scaled neighbourhood of $I_1$; \item
$F_i$ is low and central for $i=0, \dots, m$; \item there is some
$0<i<m$ with $ \frac{|I_i|}{|I_{i+1}|} <1+\delta' $,
\end{enumerate} then for $1 \le k <m$,
$$\frac{1}{C} \frac{1}{\min(k, m-k)^2} < \frac{|I_{i+k-1}\setminus
I_{i+k}|}{|I_i|} < \frac{C}{\min(k, m-k)^2}.$$ }

For similar statements see \cite{dfdm} and \cite{sha}.

\begin{proof}
We first point out the following claim.
\begin{claim}
For $f$ as in the lemma, there exists some $C(f, \delta,
\delta')>0$ \st
$$\frac{|I_m|}{|I_{0}|}> C(f, \delta, \delta').$$
\label{claim:am}
\end{claim}
This is proved in Section 5 of \cite{sha}.  One consequence of
this is that $\frac{|I_m\setminus I_{m+1}|}{|I_0|}$ is uniformly
bounded below.  This is one of the assumptions in the statement of
the Yoccoz Lemma in \cite{dfdm}.

Our proof now involves using a result of \cite{st}, the bound
$\delta$ and the small size of $I_0$, to find a nearby map in the
Epstein class. The structure of such maps, particularly at
parabolic fixed points, along with some new coordinates, give us
estimates for $\frac{|I_{i+k-1}\setminus I_{i+k}|}{|I_i|}$.

We suppose that $s>0$ is \st $F_0|_{I_1} = f^s|_{I_1}$.  We
observe that $f^{s-1}$ has uniformly bounded distortion depending
on $\delta$. We will denote $F_0|_{I_1}$ by $F$.  Letting $\psi:
[a_m,a_1] \to [0,1]$ be an affine diffeomorphism we will work with
the map $\psi\circ F\circ \psi^{-1}$.  For the rest of the
appendix we will abuse the notation and denote this map by $F$
too.

Previously we assumed that $F|_{I_1}$ had a maximum at $c$.  It
will be convenient to suppose now for this section that $c$ is a
minimum for $F|_{I_1}$.  Also we let $I_i = (a_i', a_i)$.  So in
particular, $F(a_{i+1})=a_i$.  We firstly define a point which
allows us to partition $[a_m, a_1]$ in another way.

Let $x_0 \in [a_m, a_1]$ be defined so that $|F(x_0)-x_0| =
\min_{a_m \le x \le a_0}|F(x)-x|$.  It is easy to show that
$DF(x_0)=1$.
%Note that $F(x) > x$ for all $x \in [a_m, a_1]$ and
%that $F$ is increasing.  Then let $h(x):=F(x) - (F(x_0)-x_0)$. By
%the definition of $x_0$, we have $h(x)-x>0$. Then we express $h$
%as $h(x) = x_0+(x-x_0)Dh(x_0)+ O(|x-x_0|^2)$ and so $h(x)-x =
%(x-x_0)(Dh(x_0)-1)+ O(|x-x_0|^2)$. Therefore, if $Dh(x_0)>1$ then
%there exists some $x<x_0$ near $x_0$ \st $h(x)-x<0$.  Similarly,
%if $Dh(x_0)<1$ then there exists some $x>x_0$ near $x_0$ \st
%$h(x)-x<0$.  In either case we have a contradiction, so
%$Dh(x_0)=1$.
We will suppose throughout that $|F(x_0)-x_0|$ shrinks to zero as
$|I_0| \to 0$: otherwise the proof is much simpler. We can
estimate the shape of $F$ near $x_0$ using the following
definition and lemma.

Let $\kappa>0$.  We say that the real analytic map $f:[0,1] \to
[0,1]$ is in the {\em Epstein class} $\mathcal{E}_\kappa$ if $f(x)
= \varphi Q \psi$ where $Q$ is the quadratic map $Q(z) = z^2$,
$\psi$ is an affine map and $\varphi:[0,1] \to [0,1]$ is a
diffeomorphism whose inverse has a holomorphic extension which is
univalent in the domain $\CCC_{(-\kappa, 1+\kappa)}:= \CCC
\setminus ( (-\infty,-\kappa] \cup [1+\kappa, \infty))$. For more
details on maps in this class see \cite{ms}. The following lemma
is proved in \cite{st}.

\begin{lem}
Let $f \in NF^2$.  Suppose that $I$ is a nice interval around $c$
and $J$ is a first entry domain which is disjoint from $I$ and
with entry time $s$. Suppose that $\delta>0$ is some constant \st
there exists some $\hat{J} \supset J$ \st $f^s: \hat{J} \to I'$ is
a diffeomorphism where $I'$ is a $\kappa$--scaled neighbourhood of
$I$ and $\sum |f^j(\hat{J})| \le 1$. Let $\tau_0:J \to [0,1]$ and
$\tau_s:I \to [0,1]$ be affine diffeomorphisms. Then for all
$\epsilon >0$ there exists $\delta>0$ \st $|I|< \delta$ implies
that there exists some function $G: I \to I$  in the Epstein class
$\mathcal{E}_{\frac\kappa 2}$ \st $\|\tau_s\circ
f^s\circ\tau_0^{-1} -G\|_{C^2}< \epsilon$. \label{lem:epg}
\end{lem}

We use this to prove the following claim.
\begin{claim}
There exists some $0<A<B$ \st for $I_0$ sufficiently small
$$F(x_0)+(x-x_0)+A(x-x_0)^2 \le F(x) \le
F(x_0)+(x-x_0)+B(x-x_0)^2. $$ \label{claim:AB}
\end{claim}

\begin{proof}
We know that $f^s:I_2 \to I_1$ has the following property.  The
map $f^{s-1}:f(I_2) \to I_1$ has an extension to $I_0$.
Furthermore, since $I_0$ is a $\delta$--scaled neighbourhood of
$I_1$ we use Lemma~\ref{lem:epg} to obtain a map $G_\infty$ in the
Epstein class which is $C^2$--close to $f^s$.

In fact we can choose different starting intervals $I_n$ with the
same real bounds which are smaller and smaller and which are then
rescaled to maps $F_n$ which map from the unit interval to itself.
For each such map we obtain the nearby map $G_n$ in the Epstein
class where $\|F_n - G_n\|_{C^2} \to 0$ as $n \to \infty$. For
$F_n$ we let $x_0^n$ denote a point which is equivalent to $x_0$
for $F$. Since we assume that $|F_n(x_0^n)-x_0^n|$ goes to zero,
our limit map $G_\infty$ has a parabolic fixed point $x_0^\infty$.
Also $D^2G_\infty(x_0^\infty)>0$. Thus, there exist $0<A<B$
depending only on $f$ \st for all $x \in [0,1]$ we have
$$G_\infty(x_0^\infty)+(x-x_0^\infty)+A(x-x_0^\infty)^2 \le
G_\infty(x) \le
G_\infty(x_0^\infty)+(x-x_0^\infty)+B(x-x_0^\infty)^2.$$ Clearly,
for large $n$, we have the same condition for $G_n$. Therefore, if
we take $I_0$ small enough, we may assume that it holds for $F$
too.
\end{proof}

We denote $\epsilon:= F(x_0) - x_0$. Then we have
$$\epsilon+A(x-x_0)^2 \le F(x) - x \le \epsilon+B(x-x_0)^2.$$
We suppose that $N$ is \st  $x_0 \in [a_N, a_{N+1})$.  Then for $0
\le i \le N-1$ we let $x_i := F^i(x_0)$.  We will use this
equation to find estimates for $a_j-a_{j+1}$.  Throughout we will
let $C, C'$ denote some constants depending only on $\delta,
\delta'$.

\begin{claim}
$$N =O\left(\frac{1}{\sqrt\epsilon}\right).$$ \label{claim:sqrtn}
\end{claim}

\begin{proof}
Let $N' = \max \{1 \le j \le N-1: x_j - x_0 \le \sqrt\epsilon\}$.
We will first show that $N'$ satisfies the claim.  For $j \le N'$,
we have
$$\epsilon \le x_{j+1} -x_0 \le \epsilon(B+1).$$
Therefore,
$$N' \epsilon \le \sum_{j=0}^{N'-1}x_{j+1} -x_j \le N'\epsilon(B+1).$$
Since $\sum_{j=0}^{N'-1}x_{j+1} -x_j = x_{N'}-x_0 \le
\sqrt\epsilon$ we have $N' \le \frac{1}{\sqrt\epsilon}$.
Furthermore, $x_{j+1}-x_0> \sqrt\epsilon$ so $\epsilon(N'(B+1)+1)
> \sqrt\epsilon$ and $N'> \frac{1}{(B+1)\sqrt\epsilon} -1$. I.e.
$N' =O\left(\frac{1}{\sqrt\epsilon}\right)$.

Next we find estimates for $N-N'$.  For $N'<j \le N$ we again
consider the equation
$$\epsilon+A(x_j-x_0)^2 \le x_{j+1} - x_j \le \epsilon+B(x_j-x_0)^2.$$
But note that here $B(x_j-x_0)^2> \epsilon$ so we can write
instead
$$A(x_j-x_0)^2 \le x_{j+1} - x_j \le 2B(x_j-x_0)^2.$$
We make a change of coordinates.  We let $y_j :=
\frac{1}{x_j-x_0}$.  Then we have
$$y_j - y_{j+1} = \frac{x_{j+1} -
x_j}{(x_{j}-x_0) (x_{j+1}-x_0)}. $$  By the above bounds we have
$$\frac{A(x_j-x_0)}{x_{j+1}-x_0}< y_j - y_{j+1}
< \frac{2B(x_j-x_0)}{x_{j+1}-x_0} < 2B.$$  Furthermore,
$$ y_j - y_{j+1}> \frac{A(x_j-x_0)}{(x_{j+1}-x_j)+(x_j-x_0)}>
\frac{A(x_j-x_0)}{2B(x_j-x_0)^2+(x_{j}-x_0)} > \frac{A}{2B+1}.$$
Observe that $x_N \in (a_1, a_0)$ and $|a_0-a_1|> \delta$.  So
since $|x_N-x_{N-1}|$ is approximately $|a_0-a_1|$ and since we
fixed $\delta$, we know that $y_N = O(1)$. Also note that $y_{N'}
=O\left(\frac{1}{\sqrt\epsilon}\right)$ and so $y_{N'} -y_N
=O\left(\frac{1}{\sqrt\epsilon}\right)$. Summing we obtain
$$\frac{C}{\sqrt\epsilon} < y_{N'}-y_N = \sum_{j=N-1}^{N'} y_j -
y_{j+1}< 2B(N-N')$$  and
$$\frac{C'}{\sqrt\epsilon} >  y_{N'}-y_N = \sum_{j=N-1}^{N'} y_j -
y_{j+1} >  \frac{A(N-N')}{2B+1}.$$  So
$N-N'=\left(\frac{1}{\sqrt\epsilon}\right)$ too.  Adding this to
the estimates for $N'$ we prove the claim.
\end{proof}

To prove the lemma, we will use Claims~\ref{claim:am} and
\ref{claim:sqrtn} together, along with bounded distortion, which
means that $a_j-a_{j+1}$ is like $x_{N-j} - x_{N-j-1}$.

Firstly we will use the above coordinate change again.  For $j>
N'$ we have $$y_j > y_j-y_N = \sum_{j=N-1}^{j} y_i - y_{i+1}>
\frac{A(N-j)}{2B+1}$$ and so $\frac{1}{x_j-x_0} >
\frac{A(N-j)}{2B+1}$ and $x_{j+1} - x_j <
2B\left(\frac{2B+1}{A(N-j)}\right)^2$.

We have proved that if $0 \le j \le N'$ then \begin{equation}
\epsilon < x_{j+1}-x_j <   C'\epsilon \label{eq:less than
N'}\end{equation} and if $N' < j \le N$ then \begin{equation}
\epsilon < x_{j+1}-x_j < \frac{C'}{(N-j)^2}.\label{eq:greater than
N'}\end{equation} Similarly we can define $x_j= F^j(x_0)$ for
negative $j$ where $0 \le |j|< m-N$.  Now we will show that
Claim~\ref{claim:sqrtn} follows for this situation too and we get
equivalents to \eqref{eq:less than N'} and \eqref{eq:greater than
N'}. We define some $M'$ analogously to the definition for $N'$
and so if $ |j| \le M'$ then $$ \epsilon< x_{j+1}-x_j <
C'\epsilon.$$ And if $M' < |j| \le m-N$ then
$$\frac{C}{(m-N+j)^2} < x_{j+1}-x_j < \frac{C'}{(m-N+j)^2}.$$
(In the step of the proof where estimates on $y_{N-m}$ are
required, we use Claim~\ref{claim:am} to give $|a_{m-1}-a_m|$
uniformly bounded below and the fact that $|x_{-m-1}-x_{m}|$ is
approximately $|a_{m-1}-a_m|$.)  Note also that we can show that
$m-M' =O\left(\frac{1}{\sqrt\epsilon}\right)$.

Observe that $a_j-a_{j+1}$ is essentially the same as $x_{N-j} -
x_{N-j-1}$. So if $N \ge j\ge N-N'$, we have
$$C\epsilon < a_j-a_{j+1} <   C'\epsilon.$$  Observe that
$\frac{1}{N-N'} \ge \frac{1}{j} \ge \frac{1}{N}$.  Since $\epsilon
\left(\frac{1}{N^2}\right)$ and $\epsilon
\left(\frac{1}{(N-N')^2}\right)$ this implies that we have
$$\frac{C}{j^2}< a_j-a_{j+1} < \frac{C'}{j^2}.$$

Now  if $N-N' \ge  j \ge O(1)$ then clearly we have $a_j-a_{j+1}<
\frac{C'}{j^2}$.  Also, \begin{align*} x_{N-j}-x_{N-j-1} &
>A(x_{N-j-1}-x_0)^2 = A\left(\sum_{k=j-1}^{N-1}
x_{N-k}-x_{N-k-1}\right)^2\\
& \ge A\left(\sum_{k=1}^{N'} x_k-x_{k-1}\right)^2 \ge A(N'\sqrt
\epsilon)^2.\end{align*} Now since
$\sqrt\epsilon=O\left(\frac{1}{N'}\right)$, we have
$x_{N-j}-x_{N-j-1} \gtrsim 1.$ Thus
$$\frac{C}{j^2} < a_j-a_{j+1}< \frac{C'}{j^2}.$$
If $ N \le  j \le m- M'$ then again we have
$$C\epsilon< a_j- a_{j+1} < C'\epsilon.$$

Note that we also have $m-N \ge m-j \ge m-M'$.  Since $m-N, m-M'
=\left(\frac{1}{\sqrt\epsilon}\right)$ we have
$$\frac{C}{(m-j)^2}< a_j - a_{j+1} < \frac{C'}{(m-j)^2}.$$
If $ m-M' \le  j \le m-1$ we have $$\frac{C}{(m-j)^2}< a_j -
a_{j+1} < \frac{C'}{(m-j)^2}$$ where the lower bound follows as
above.

To conclude, if $1 \le j \le N$ then we have some constant $C$ \st
$j \le C(m-j)$ and  $a_j - a_{j+1} \asymp \frac{1}{j^2}$.  If $N
\le j \le m-1$ then we have some constant $C'$ \st $m-j \le C'j$
and $a_j - a_{j+1} \asymp \frac{1}{(m-j)^2}$.  So in either case
we have
$$ a_j - a_{j+1} \asymp \frac{1}{(\min(j, m-j))^2}$$ as required.
\end{proof}

\medskip
\noindent
Department of Mathematics\\
University of Surrey\\
Guildford, Surrey, GU2 7XH\\
UK\\
\texttt{m.todd@surrey.ac.uk}\\
\texttt{http://www.maths.surrey.ac.uk/showstaff?M.Todd}

\medskip

\end{document}